\documentclass[12pt,leqno,openany]{book}

\usepackage{makeidx}
\makeindex

\def\diam{\mathop{\rm diam}}
\def\dist{\mathop{\rm dist}}

\newtheorem{theorem}{Theorem}
\newtheorem{lemma}[theorem]{Lemma}
\newtheorem{proposition}[theorem]{Proposition}
\newtheorem{sublemma}[theorem]{Sublemma}
\newtheorem{definition}[theorem]{Definition}
\newtheorem{corollary}[theorem]{Corollary}
\newtheorem{problem}[theorem]{Problem}
\newtheorem{remark}[theorem]{Remark}
\newtheorem{claim}[theorem]{Claim}
\newtheorem{assumptions}[theorem]{Assumptions}
\newtheorem{examples}[theorem]{Examples}
\newtheorem{question}[theorem]{Question}

\newcommand{\begintheorem}{\addtocounter{equation}{1}\begin{theorem}}
\newcommand{\beginlemma}{\addtocounter{equation}{1}\begin{lemma}}
\newcommand{\beginproposition}{\addtocounter{equation}{1}\begin{proposition}}
\newcommand{\beginsublemma}{\addtocounter{equation}{1}\begin{sublemma}}
\newcommand{\begindefinition}{\addtocounter{equation}{1}\begin{definition}}
\newcommand{\begincorollary}{\addtocounter{equation}{1}\begin{corollary}}
\newcommand{\beginproblem}{\addtocounter{equation}{1}\begin{problem}}
\newcommand{\beginremark}{\addtocounter{equation}{1}\begin{remark}}
\newcommand{\beginclaim}{\addtocounter{equation}{1}\begin{claim}}
\newcommand{\beginassumptions}{\addtocounter{equation}{1}\begin{assumptions}}
\newcommand{\beginexamples}{\addtocounter{equation}{1}\begin{examples}}
\newcommand{\beginquestion}{\addtocounter{equation}{1}\begin{question}}

\begin{document}

\frontmatter

\title{Where the Buffalo Roam: \\ Infinite Processes and Infinite
Complexity}

\author{Stephen Semmes  \\
	Department of Mathematics \\
	Rice University}
\date{}

\maketitle

\chapter{Preface}


\renewcommand{\thefootnote}{}   

\footnotetext{This manuscript (modulo a few modest changes) was
distributed informally as IHES preprint M/96/77.  The author was
partially supported by the National Science Foundation.}

	Much attention has been paid in recent years to questions of
{\em complexity}, how much time or memory does it take to make a
certain construction, or to answer a particular question, or what is
the length of the shortest path with prescribed properties, and so
forth.  In mathematics one allows infinite processes and infinite
constructions, and the time seems to be ripe to look at measurements
of their complexity.  

	We shall try to paint a picture through examples.

\tableofcontents

\mainmatter

\chapter{Introduction}
\label{Introduction}

	In the formalization of mathematics there are obvious
philosophical questions of what is real and what is mere abstraction,
and how mathematics might treat the two differently.  In the daily
life of the working mathematician there is not much point to the
broader debate; the matter is typically clear in actual practice.  We
know when we have given an argument which is nonconstructive, and
sometimes we regret it.

	Still it is interesting to reflect on the general question.
We collect here some examples.  In particular we shall be interested
in asking whether or not a given mathematical statement or concept or
proof has a natural finite version at all.  Sometimes this is not
completely clear, and may suggest something new.

	Another interesting point concerns the ``complexity'' of
proofs.  In classical predicate logic one knows that the size of a
proof may have to be very large compared to the size of the statement.
In the context of propositional logic this matter is connected to the
problem ``NP = co-NP?'' \index{NP=co-NP?} in complexity theory
\cite{CR}.  (The survey \cite{CS1} provides additional information and
references.)  One might wonder about the need for infinite processes
in a proof which are substantially more complicated than the infinite
processes required for the statement.  In the context of arithmetic
one sometimes measures the infinite complexity of a proof in terms of
the strength of the induction that is used.  One can also see the
matter at play in analysis, geometry, and topology.

	Another point about infinite processes is that they lead to
more ``symmetry''.  This symmetry can take a directly geometric form,
as in the discussion of the proof of the John-Nirenberg theorem and
its finite versions in \cite{CS1}, or it can be a symmetry in
language.  One might ask only whether some quantities are finite or
infinite, for instance, and not worry about explicit bounds.  This
simplification in the language can make it easier to use and formulate
``lemmas'', and to find short proofs of theorems through the repeated
application of lemmas.  This is part of the matter of {\em cuts} and
{\em cut elimination} \index{cut-elimination} in logic, as described
in \cite{CS1}.  

	These considerations of logic, complexity, and the lengths of
proofs provided part of the motivation for looking more carefully at
infinite processes and infinite complexity, but there were other
motivations coming from ordinary mathematical activity.  It turns out
that infinite levels of complexity arise naturally and in some sense
unavoidably in geometry and topology, for instance.  We shall see more
about this later.  An enormous amount of classical analysis
is about controlling infinite levels of complexity as well.

	 Of course there are many well-known examples from number
theory and combinatorics related to infinite complexity, but we shall
emphasize geometry, topology, and analysis here.

	It is also pleasant to look for concrete examples from
ordinary mathematics to illustrate general phenomena and ideas from
formal logic.  See \cite{CS3} for further discussions of this nature,
with connections with logic more closely drawn.

	See \cite{CS2, CS3, CS4} for more discussion about the
relationship between short proofs of large objects and the presence of
symmetry.  We shall return to some of these themes in Section
\ref{Comments}.

	These questions about infinite processes and infinite
complexity provide a common thread by which to take a look at
mathematics in a broad way.  By the end we shall say a little about
many different topics, rather than the other way around.

\chapter{Examples from the infinite world}
\label{Examples from the infinite world}

	Let us look at some concrete examples from ordinary mathematics
where infinite processes are involved, perhaps unavoidably.

\section{Continuity}
\label{Continuity}
\index{continuity}

	Consider the notion of continuity, for real-valued functions
on the unit interval $[0,1]$, say.  A function $f : [0,1] \to {\bf R}$
is {\em continuous} if for every $\epsilon > 0$ and every $x \in
[0,1]$ there exists $\delta > 0$ so that
\begin{equation}
	|f(x) - f(y)| < \epsilon \quad\hbox{whenever}\quad |x-y| < \delta.
\end{equation}

	This notion seems to be unavoidably an infinite one.  It makes
no sense in a finite approximation.

	This changes if one prescribes a modulus of continuity,
\index{modulus of continuity} which is to say a particular rate at
which the effects of continuity occur.  One might consider the
Lipschitz condition \index{Lipschitz condition} for a function $f :
[0,1]
\to {\bf R}$, which asks that there be a constant $C > 0$ such that
\begin{equation}
\label{Lipschitz}
	|f(x) - f(y)| \le C \, |x-y| \qquad\hbox{for all } x, y \in [0,1].
\end{equation}
One might instead choose an $\alpha \in (0,1)$ and ask for the weaker
condition that there exist a constant $C$ such that \index{H\"older continuity}
\begin{equation}
\label{Holder}
	|f(x) - f(y)| \le C \, |x-y|^\alpha
				 \qquad\hbox{for all } x, y \in [0,1].
\end{equation}
These provide different ways to measure the continuity of a function.
They are not exhaustive; one can give examples of continuous functions
which do not satisfy any of these conditions, or any like them.

	If we do not specify a rate then it is difficult to see how we
can think of continuity intelligently in a finite context.  If we
choose a particular rate then we can do it.  Fix a large integer $N$,
and let us think of the set $I_N = \{ \frac{j}{N} : j = 0, \ldots,
N\}$ as being a finite approximation to $I = [0,1]$.  Given a function
$f : I_N \to {\bf R}$, we can look at the smallest constant $C$ for
which the analogue of (\ref{Lipschitz}) is satisfied.  This is called
the {\em Lipschitz norm} of $f$.  Alternatively we can fix an exponent
$\alpha \in (0,1)$ and look at the smallest constant $C$ for which the
analogue of (\ref{Holder}) is satisfied.  We could fix other rates and
measure the approximate continuity for them.  As $N \to \infty$ these
measurements can become dramatically different.

	For each fixed $N$ no matter what rate we use we are always
working with the same functions on $I_N$, namely all functions.  When
we pass to the ``infinite'' realm of functions on $I$ these classes
become different.  For instance the function $f(x) = \sqrt{x}$
satisfies (\ref{Holder}) when $\alpha \le \frac{1}{2}$ but not when
$\alpha > \frac{1}{2}$.  

	There are amusing points here about quantifiers.  In the
infinite context of the actual unit interval $[0,1]$ it makes sense to
talk about the general notion of continuity, to allow all rates of
continuity at once.  This universality is accommodated through the use
of {\em universal} and {\em existential quantifiers}, in the sense of
formal logic, and does not work so well for an individual $I_N$.

	In general the infinite context seems to provide a setting
where one can work with quantifiers more efficiently.  This point
interacts with the idea of {\em cut elimination} from formal logic in
any interesting way.  (See \cite{CS1} for an introduction.)  One knows
that quantifiers and their nesting can be conducive to large expansion
in cut elimination.  One might say instead that the use of quantifiers
can help one to make short proofs of complicated statements.

\section{Uniform continuity}

	Remember that continuous functions on the unit interval are
always uniformly continuous.  \index{uniform continuity} That is,
continuity implies that for each $\epsilon > 0$ there is a $\delta >
0$ so that
\begin{equation}
	|f(x) - f(y)| < \epsilon \quad\hbox{ whenever }\quad |x-y| < \delta,
\end{equation}
i.e., there is always a uniform rate at which continuity takes effect.
Does this theorem have a natural ``finite'' version?  It would seem
not.  This is reasonable, since we need to be in the infinite context
to make sense of either the assumption or the conclusion.  In other
words this is a theorem about the relationship between two infinite
concepts, with perhaps nothing to say about finite mathematics.

	This is a bit troubling.

	If we fix a particular (uniform) rate of continuity, as in the
Lipschitz condition (\ref{Lipschitz}), or the H\"older condition
(\ref{Holder}), then we do not really lose information in passing
between the infinite context and its finite approximations.  For
instance, given a Lipschitz function on $I$, its restriction to each
$I_N$ is Lipschitz with the same constant; for large $N$ this
restriction determines the original function up to a small error that
one can easily compute; any function on $I_N$ can be extended to give
a function on $I$ with the same Lipschitz constant, etc.  However,
having a uniform rate of continuity to begin with is very different
from the ordinary notion, which allows different rates at different
points.

	This leads to another interesting point about continuity: How
does one check it in practice?  Not in a circular situation, like when
one assumes the continuity of two functions and wants to prove the
continuity of the sum, but more seriously where an actual function is
under consideration.  I do not know of a concrete construction which
implies the continuity of a function without the possibility to obtain
information about the actual rate of continuity that one can somehow
write down, at least in principle (and perhaps not easily).  The word
``concrete'' might necessarily exclude many constructions from
functional analysis, such as the Hahn-Banach theorem and compactness
arguments, as well as transversality arguments and the Baire category
theorem.  We shall return to these methods later in the text.

	In my experience even some not-so-concrete situations lead to
particular degrees of continuity.  I am thinking for instance of
solutions of partial differential equations which are obtained by
minimizing some kind of ``energy'' functional.  This means that the
function may have been obtained nonconstructively through a
compactness argument.  Still one can often prove that such a function
is continuous, or that its first derivatives are continuous, at least
on the complement of a small ``singular'' set.  The methods used in
the field that I know do not give continuity with no other
information, they tend to provide a condition like (\ref{Holder}) (at
least locally) or they do not work at all.  See \cite{Gi}, for
instance.  In many important cases continuity implies smoothness, and
the issue degenerates in the end even if the idea makes sense at an
earlier stage of the argument.  This does not always happen though,
particularly in the more complicated situation where one is trying to
study the structure (like smoothness) of a singular set.

	Another interesting aspect of continuity is provided by the
theorem which says that a continuous function on a compact space
always assumes a maximum and a minimum.  What is the ``finite''
version of this?  The fact that a function on a finite set always
assumes a maximum and a minimum is not very interesting.  In some
sense continuity and compactness exactly provide the tools to reduce
the ``infinite'' problem to a finite one, and the finite and infinite
versions are not really so different.

	Alternatively one might think just about the fact that a
continuous function on a compact space is automatically bounded, or
that if it is positive everywhere then it is actually bounded away
from $0$.  In other words, the existence of uniform bounds which are
trivial for finite sets but not in infinite contexts.  

	Note well the difference between ``uniform'' and ``universal''
bounds.  A real-valued function on a finite set is always uniformly
bounded, but there is no universal bound.  The uniform boundedness of
continuous functions on compact sets straddles the two.  

	Sometimes one derives universal bounds in concrete situations
through a compactness argument of this type, with no clue as to
explicit constants.  In practice these arguments typically proceed by
contradiction.  One assumes that the existence of a sequence for which
there is no bound, one uses compactness to extract a subsequence which
converges to something, one uses information in the problem to say
something about the limit.  Effectively one is proving some kind of
continuity, but sometimes one can get away with a bit less than that.
Arguments of this type tend to be a bit mysterious.

\section{Classical Banach spaces}

	A {\em Banach space} \index{Banach space} is a normed vector
space (over the real or complex numbers) which is {\em complete} as a
metric space (using the metric induced by the norm).  Remember that a
{\em norm} on a vector space $V$ is a nonnegative real-valued function
$\|\cdot\|$ on $V$ which vanishes only at the zero vector and satisfies
the property that
satisfies 
\begin{equation}
	\| \lambda \, v \| = |\lambda| \, \|v\|
\end{equation}
for all vectors $v$ and scalars $\lambda$, and also the {\em triangle
inequality}
\begin{equation}
	\|v+w\| \le \|v\| + \|w\|
\end{equation}
for all vectors $v$ and $w$ in $V$.  The {\em metric induced by the norm}
is defined by $d(v,w) = \|v-w\|$.

	Basic examples include the $\ell^p$ spaces \index{$\ell^p$
spaces}, $1 \le p \le \infty$, which are the spaces of sequences $x =
\{x_j\}_{j=1}^\infty$ such that
\begin{equation}
	\|x\|_p = (\sum_j |x_j|^p)^\frac{1}{p} < \infty.
\end{equation}
When $p = \infty$ we take $\ell^\infty$ to be the space of bounded sequences,
with norm
\begin{equation}
	\|x\|_\infty = \sup_j |x_j|.
\end{equation}
It is well known that these are all norms and that the associated
metric spaces are complete.

	For $p=\infty$ there is another natural space, namely $c_0$,
\index{$c_0$} which is the space of sequences $x =
\{x_j\}_{j=1}^\infty$ such that
\begin{equation}
	\lim_{j \to \infty} x_j = 0.
\end{equation}
For this space we use the norm $\| \cdot \|_\infty$ again, and one can
still show completeness.  

	Notice that the space of sequences which have zero for all but
finitely many coordinates is dense in $\ell^p$ for $1 \le p <
\infty$.  This is not true when $p = \infty$.  The closure in $\ell^\infty$
of this space is precisely $c_0$.

	Given a Banach space $(X, \|\cdot\|)$ we define its {\em dual}
\index{dual of a Banach space}
$(X^*, \|\cdot\|_*)$ to be the space of bounded linear functionals on
$X$.  A {\em linear functional} simply means a linear mapping to the
ground field (${\bf R}$ or ${\bf C}$ in this case).  A linear
functional $\lambda$ on $X$ is said to be {\em bounded} \index{bounded
linear functionals} if there is a constant $M < \infty$ such that
\begin{equation}
	|\lambda(x)| \le M \, \|x\| 	\qquad\hbox{for all } x \in X.
\end{equation}
The smallest such $M$ can be written as
\begin{equation}
	\|\lambda\|_* : = \sup_{x \in X \atop \|x\| = 1} |\lambda(x)|
\end{equation}
and defines a norm on the space of bounded linear functionals.  With
this norm $(X^*, \|\cdot\|_*)$ becomes a Banach space in its own
right (i.e., it is complete).

	If $1 \le p < \infty$ then the dual of $\ell^p$ is $\ell^q$,
where $q$ is determined by the equation $\frac{1}{p} + \frac{1}{q} =
1$.  (Thus $q = \infty$ when $p = 1$.)  In particular if $1 < p <
\infty$ then the dual of the dual of $\ell^p$ is $\ell^p$ again.
This breaks down for $p = 1, \infty$.  The dual of $\ell^\infty$ is a
mess and it is better not to talk about it.  However the dual of $c_0$
is $\ell^1$.  Thus the second dual of $c_0$ is $\ell^\infty$, a nice
feature which is rather basic in analysis.

	What about {\em finite-dimensional} spaces?  For them
one can show that the second dual of a space is always the original space
again.

	Fix an integer $n$, and define $\ell^p_n$ in the same way as
above, except that we restrict ourselves now to sequences of length
$n$.  Thus all of our spaces are just ${\bf R}^n$ or ${\bf C}^n$ as
vector spaces, the differences are in the norms, which we define
exactly as above.  We cannot define $c_0$ now, or rather we cannot
make a distinction between $c_0$ and $\ell^\infty$; there is no limit
to take.  This fits well with the story of duality, because each of
$\ell^1_n$ and $\ell^\infty_n$ is the Banach-space dual of the other.
The content in this statement lies in the norms, since the underlying
vector spaces and their duals are just ${\bf R}^n$ or ${\bf C}^n$.

	Thus again we have some distinctions in the ``infinite'' case
which do not make sense in the ``finite'' approximations.

	Although in finite dimensions the second dual of a space always
gives the space back, one can make quantitative distinctions between
the cases of $p=1, \infty$ and $1 < p < \infty$ for the $\ell^p_n$'s.
For instance, the norms for $1 < p < \infty$ enjoy stronger convexity
properties which persist to infinite dimensions.  

	See \cite{HS, Ru2, Ru3} for more information about these topics
(at an introductory level).

\section{The open mapping theorem}

\index{open mapping theorem}

\begin{theorem}
Let $(X, \|\cdot\|_X)$ and $(Y, \|\cdot\|_Y)$ be Banach spaces, and suppose
that $T : X \to Y$ is a bounded linear mapping.  (Boundedness means that
there is a constant $M > 0$ so that
\begin{equation}
	\|T(x)\|_Y \le M \, \|x\|_X \qquad\hbox{for all } x \in X,
\end{equation}
and this condition is equivalent to continuity in the case of {\em
linear} mappings.)  If the image of $T$ is all of $Y$, then $T$ is an
open mapping, i.e., it sends open subsets of $X$ to open subsets of
$Y$.
\end{theorem}

	This is a well-known theorem due to Banach.  (See \cite{Ru2}.)
The main point about the assumption of surjectivity is that the image
of $T$ be closed.  If it is closed then one can always treat the image
as a Banach space in its own right.

	To understand the meaning of this theorem it is helpful to
consider some special situations.  The first point is that the theorem
does not have much content for finite-dimensional spaces.  In finite
dimensions all norms give equivalent topologies, and the theorem
reduces to a simple observation about linear mappings between
Euclidean spaces.  In infinite dimensions it is much more substantial.

	What does the condition that $T$ be open really mean?  Using
linearity one can show that $T$ is an open mapping if and only if
$T(B_X(0,1))$ contains the origin of $Y$ in its interior, where
$B_X(0,1)$ denotes the unit ball in $X$.  In other words, one only has
to test for openness ``once''.  This is the same as asking that there
exist an $\eta > 0$ so that
\begin{equation}
\label{open at 0}
	T(B_X(0,1)) \supseteq B_Y(0,\eta).
\end{equation}
Assume now that $T$ is one-to-one.  In this case (\ref{open at 0}) is
equivalent to the lower bound
\begin{equation}
\label{lower bound}
	\|T(x)\|_Y \ge \eta \, \|x\|_X \qquad\hbox{for all } x \in X.
\end{equation}

	This is very nice.  The lower bound (\ref{lower bound})
implies that the kernel of $T$ is trivial, and it is a kind of
quantitative version of injectivity.  It is an automatic consequence
of injectivity in finite dimensions but not in infinite dimensions,
and the content of the open mapping theorem is that in infinite
dimensions (\ref{lower bound}) follows from injectivity if one also
knows surjectivity.  

	As noted above surjectivity can be replaced by the requirement
that the range of $T$ be closed, and indeed (\ref{lower bound}) forces
the range to be closed.  The latter assertion is a standard exercise;
the (metric) completeness of the domain implies the completeness of
the image, under the assumption of a lower bound like (\ref{lower
bound}).  Thus one can also think of (\ref{lower bound}) as being a
quantitative version of the requirement that the range be closed.

	The open mapping theorem can therefore be viewed as a
criterion for the existence of a uniform bound.  As such it is a bit
like the fact that continuous functions on compact spaces are always
uniformly bounded.  

	To make the role of the uniform bound more concrete let us
consider a special case.  Take $X = Y = \ell^p$ as defined in the
previous section, with any (single) choice of $1 \le p \le
\infty$.  We could also take $X = Y = c_0$ instead.  Fix a sequence of
scalars $\{t_j\}_{j=1}^\infty$ which are bounded and never $0$.
Define an operator $T$ by
\begin{equation}
	T(x) = y, \enspace y_j = t_j \, x_j \quad\hbox{for all } \, j.
\end{equation}
That is, $x$ and $y$ are sequences, and $T$ is a ``diagonal'' operator
with entries $t_j$.  This defines a bounded linear operator on $X = Y
= (\ell^p$ or $c_0)$.  The requirement that the $t_j$'s never vanish
is the same as the injectivity of the operator $T$, and it is also
equivalent to the requirement that the range of $T$ be dense in this
special case of diagonal operators.  One can check that (\ref{lower
bound}) is the same as
\begin{equation}
\label{lower bound for scalars}
	|t_j| \ge \eta \qquad\hbox{for all } \, j
\end{equation}
in this case.  Thus the open mapping theorem provides a criterion
for a uniform lower bound on the size of these scalars $t_j$.

	If we were to make the analogous construction in finite
dimensions, using scalars $\{t_j\}_{j=1}^n$ and defining $T$ on
$\ell^p_n$, we would simply say that (\ref{lower bound for scalars})
holds as soon as we know that the $t_j$'s are all different from $0$.
In finite dimensions surjectivity is the same as having dense image
for linear mappings.  In infinite dimensions these properties are more
subtle and involve infinite processes.  The effect of the open mapping
theorem is to relate two different kinds of infinite processes in a
nontrivial way.

	Is the open mapping theorem useful?  It is frequently very
important in analysis to be able to derive lower bounds like
(\ref{lower bound}).  In my experience however the open mapping
theorem has been primarily of psychological value.  The main approach
that an analyst has to proving something like the surjectivity of a
mapping is to first establish an inequality like (\ref{lower bound}).
The open mapping theorem implies in some sense that there is no other
way to do it.

\section{Category versus measure}

	In mathematics it is very convenient to have reasonable ways
to talk about ``almost all'' points in a space.  There are two
especially famous ways to do this.  The first is to use the Baire
category theorem. \index{Baire category theorem} In a complete metric
space we might think of a set $E$ which is a countable intersection of
dense open sets as containing almost all elements of the original
space.  Baire's theorem implies that any such $E$ is always dense.  In
particular the intersection of any two or even countably many of these
putatively large sets remains dense and in the same class of subsets.

	Another possibility is to work on a measure space and to look
at sets whose complement has measure zero.  Again we can take {\em
countable} intersections and still have a set of the same type.

	In some cases one can use algebraic or analytic structure to
make much more refined notions of ``almost all'' points in a space.
One can use proper subvarieties to make notions of negligibility for
subsets.

	A more primitive idea is to use cardinality, to simply look at
subsets whose complement has cardinality bounded in a certain way.
Inside an infinite set one might consider subsets whose complement is
finite.  For uncountable sets one can look at subsets whose complement
is at most countable.

	In metric spaces one can impose more precise ``smallness''
conditions on the complement of a set, through Hausdorff measure or
Minkowski content, for instance.  In topology one might use
submanifolds of smaller dimension, or just subsets of smaller
topological dimension.  (These concepts will be defined later in the
text.)

	The idea of ``almost every'' cannot work in the finite world
in a similarly nice way.  We cannot have notions of ``almost every''
which are closed under arbitrary finite intersections and which are
also nontrivial.  The empty set should not contain almost all the
elements of a nonempty set!  What can we do in finite contexts?

	One possibility is to replace ``zero measure'' with ``small
measure'', and then restrict the kind of intersections that one makes.
This is awkward but feasible nonetheless, and commonplace in counting
arguments.  The Baire category theorem is more problematic.  I do not
see a natural version for finite mathematics.  This is slightly
disturbing, particularly for the idea of using Baire category as a
method for proving existence of points with special properties.

	The proof is less scary.  One ``finds'' actual sequences which
converge to a point in the intersection.  The method is fairly
concrete even if it is not so clear how to make more concrete
statements.  To analyze it carefully one can think about trying to
implement it on a computer, to actually produce the relevant sequence.
We then have to implement our ``assumptions'' that the initial sets
are dense and open.

	This is not ridiculous, just very awkward.  Here again we see
how quantifiers can be manipulated much more efficiently in
``infinite'' contexts, and how much more simple and ``symmetric''
definitions can be in the infinite world.  The notions of openness and
density rely on quantifiers very explicitly.  One also sees some of
the cycling associated to proofs with cuts and quantifiers.

	For the record let us sketch the proof.  We have a complete
metric space $(M, d(x,y))$ and a sequence $\{U_j\}_{j=1}^\infty$
of dense open sets, and we want to show that $\bigcap_j U_j$
is also dense.

	A basic lemma is that $V_j = U_1 \cap U_2 \cap \ldots \cap
U_j$ is dense and open for each $j$.  One can check that the
intersection of a pair of dense open sets is dense and open, and then
repeat.  Although simple to prove this part is already substantial, in
the sense that dense open sets provide a cruder notion of ``almost
all'' points in a metric space with good properties.  This cruder
notion also lacks an obvious counterpart for finite mathematics,
and enjoys a kind of symmetry which is useful in making proofs.

	Next one chooses a sequence of points $\{x_j\}_{j=1}^\infty$
in $M$ and a sequence of positive numbers
$\{\epsilon_j\}_{j=1}^\infty$ with the following properties:
\begin{enumerate}
\item		each $x_j$ lies in $V_j$, and in fact 
		$\overline{B}(x_j, 2 \, \epsilon_j) \subseteq V_j$,

\item		$d(x_{j+1}. x_j) \le \epsilon_j$ for all $j$, and

\item		$\epsilon_{j+1} \le \frac{\epsilon_j}{2}$ for all $j$.

\end{enumerate}
It is not difficult to make these choices, using the fact that
$\{V_j\}_{j=1}^\infty$ is a decreasing sequence of dense open subsets.
They imply that 
\[
	d(x_j, x_k) \le 2 \, \epsilon_j \quad\hbox{ when } k \ge j,
\]
by summing a geometric series.  Thus $\{x_j\}_{j=1}^\infty$ is a
Cauchy sequence in $M$, and therefore converges, and the properties
above ensure that the limit lies in all the $V_j$'s.

	Note that the notions of ``almost every'' provided by Baire
category and measure need not be compatible with each other, even on
the real line using Lebesgue measure.  This is true in principle and
in actual practice.

	If one wanted to analyze Baire category arguments more
carefully, looking for more quantitative versions, or just to bring
out some kind of constructivity, then a natural class of spaces to
consider would be the following.  Let $\{F_j\}$ be some sequence of
{\em finite} sets, and let ${\cal F}$ denote their Cartesian product.
Thus an element $x \in {\cal F}$ actually represents a sequence $\{x_j\}$,
for which $x_j \in F_j$ for each $j$.  This space can be given a topology
in the usual way, where two elements $x$ and $y$ are ``close'' if their
entries $x_j, y_j \in F_j$ agree for a large number of $j$'s.  This is
the usual {\em product topology}, and in this case it defines a space
which is compact and for which the topology is determined by a metric.

	Now in this case the properties of being dense or open admit more
concrete interpretations than usual, and one can analyze the argument above
to make it more concrete.

\chapter{Symmetry and language}
\label{About symmetry in language}

	There is a common thread to the examples above, concerning
``abstract'' concepts in infinite contexts and the possibility for
meaningful ``concrete'' versions in finite situations.  Here is a
sampling.

\begin{itemize}

\item 	We have {\em continuity} and {\em uniform continuity} for functions,
	and then more concrete versions where a {\em rate} of continuity
	is specified, as in the Lipschitz and H\"older conditions.

\item  	We can talk about the {\em boundedness} of linear mappings 
	between Banach spaces, and then more concretely we can look at
	the {\em size} of the bounds.  In particular {\em all} linear
	operators between finite-dimensional spaces are bounded.

\item 	We can think of a subset of a space as being negligible
	if it has measure zero, while in the finite world we might
	have to be more careful and consider the actual number of
	elements of the subset as compared to the ambient space as a
	whole.

\item	In topological spaces we can look at the properties of being 
	{\em dense} and {\em open}, for which elegant concrete versions
	are not apparent.

\item	We can consider merely the {\em finiteness} or {\em infiniteness}
	of some quantity, as opposed to trying to decide which numbers
	are {\em large} and which are {\em small}.

\end{itemize}

\noindent
This list is by no means complete, but it gives a flavor.  The extent
to which the infinite concepts can be captured by finite versions
depends on the details of the situtation.

	The extent of abstractness or concreteness used has important
concequences for the enterprise of making proofs.  In mathematics we
consider something to be ``known'' when we have proved it and not
otherwise.  What if something is true but all proofs are very long?

	Mathematicians have long been concerned about matters of
undecidability.  We have the goal of ``truth'', we try to formalize
it, we make axiom systems to cover what we are interested in, and then
some questions are typically left over.  This is not too surprising
after all, but we often enter the field as students with more idealism
and na\"{\i}vet\'e about the notion of ``truth''.  There are even
precise theorems of G\"odel and others to the effect that this is
unavoidable in the context of arithmetic.  (See
\cite{Ma}, for instance.)  There are even quantitative versions of 
G\"odel's incompleteness theorems which have been given by Chaitin,
based on algorithmic information theory.

	In simpler forms of logic we do not have this problem.
Propositional and first-order predicate logic are consistent and
complete.  ``All is well'', in a sense.

	But is it really?  What about ``practical'' incompleteness,
where a statement is true but all the proofs are too enormous for a
human being to find?

	In predicate logic it is well known that the set of
tautologies form an {\em algorithmically undecidable set}.  This means
that there is no computer program which can input a sentence in
(first-order) predicate logic and decide whether or not the sentence
is a tautology.  There are computer programs that can generate a
complete list of all (first-order) tautologies; one can simply list
all strings that could possibly represent a proof, check which ones do
represent a proof, keep those and throw away the rest.  The problem is
to know whether a given sentence will appear on the list eventually or
not, and for that there is no such algorithm.  See \cite{EFT}, for
instance.

	As a practical matter this means that there can be tautologies
whose proofs are enormously larger in size than the size of the
statement itself.  ``Enormously'' means here larger than any recursive
function in the size of the statement of the tautology.

	There are similar issues for propositional logic (i.e.,
classical logic but without quantifiers, relations, etc.).  In that
case the set of tautologies is decidable and quite easily so, one has
only to check truth tables.  The problem is that the na\"{\i}ve
algorithm could take an exponential amount of time (in the size of the
given propositional formula) to decide whether the formula is a
tautology or not.  Exponential expansion is mild compared to
nonrecursive growth, but still prohibitive practically.

	There is an algorithm which can determine in polynomial time
whether a given propositional formula is a tautology if and only if P
= NP.  The latter is a famous question in complexity theory; see
\cite{GJ, HU, Pap} for background information about it and
related problems of algorithmic complexity.  The more precise
statement is that the set of propositional formulas which are {\em
not} tautologies form an NP-complete set.  As mentioned in Chapter
\ref{Introduction}, there is a proof system in which propositional
tautologies always have polynomial-sized proofs if and only if NP =
co-NP; see \cite{CR} for more details.  

	Even in propositional logic one has therefore serious
difficulty with being able to find proofs of modest size for all
tautologies.

	What is it that allows one to make short proofs?  In
\cite{CS2, CS3, CS4} there is an idea put forward
to the effect that the existence of short proofs of complicated
statements reflects the presence of some kind of underlying symmetry.
This might be a literal ``geometric'' symmetry of some object, but it
might also appear as a symmetry in the language, in the way that
objects are described.  In our exploration of infinite processes and
infinite concepts in mathematics we see many examples which illustrate
this idea from \cite{CS2, CS3, CS4}.

	It is easier to make reasoning about the ``abstract'' concept
of finiteness or infiniteness of some quantity than it is to make more
precise estimates.  It is easier to formulate lemmas based on the
distinction between finiteness and infiniteness, and to use the same
lemma many times.  

	These are familiar themes to any working mathematician.  One
knows from practical experience that it is easier to work with the
notion of {\em continuity} of functions that it is to keep track of
how the modulus of continuity changes under various operations, and
that it is easier to say that a set is small because it is the
countable union of sets of measure zero than to make some concrete
measurement of its negligibility.

	Thus it can be convenient to work with infinite processes or
concepts even if in the end we are really interested in something more
concrete.  The relevant distinctions (between large and small, for
instance) can simplify in the infinite world; we might ask for the
mere {\em existence} of a bound, without trying to keep track of what
it is.  This makes it easier to formulate and use lemmas and then to
find actual proofs.

	After we find a proof we can try to go back and analyze it to
get explicit bounds, or at least to say that such bounds are implicit
in the argument.  This idea is made general and precise in formal
logic through {\em cut elimination}, \index{cut elimination} which in
effect seeks to transform a given proof into a ``direct'' one, with no
lemmas, each step being explicit for the particular case at hand.  See
\cite{CS1} for an introduction and \cite{Gir} for more details.

\chapter{The existence of limits}

	There are many results in classical analysis and related areas
of mathematics which assert the existence of a limit, perhaps only
``almost everywhere''.  To what extent do these statements have
natural finite versions?  Does their content lie purely within the
infinite world?

\section{The Riemann-Lebesgue lemma}
\label{The Riemann-Lebesgue lemma}
\index{Riemann-Lebesgue lemma}

	Let $f(x)$ be an integrable function on the real line.  Define
its {\em Fourier transform} \index{Fourier transform} $\widehat{f}$ by
\begin{equation}
	\widehat{f}(\xi) = 
		\int_{\bf R} f(x) \, e^{- 2 \pi i \, x \, \xi} \, dx.
\end{equation}
Roughly speaking, the Fourier transform tells us how to decompose a
function into primitive ``frequencies''.  This idea is made precise
by the ``inversion'' formula
\begin{equation}
	f(x) = \int_{\bf R} 
			\widehat{f}(\xi) \, e^{2 \pi i \, x \, \xi} \, d\xi.
\end{equation}
This formula works when $\widehat{f}$ is also integrable, but in general
it requires some interpretation.  See \cite{HS, Ru2, Ru3, SW} for more
information.

	The following is a basic fact in Fourier analysis.

\begin{theorem} For each integrable function $f$ on ${\bf R}$ we have
that
\begin{equation}
	\lim_{|\xi| \to \infty} \widehat{f}(\xi) = 0.
\end{equation}
\end{theorem}

	An annoying feature of this theorem is that we cannot say
anything about how fast the Fourier transform decays at infinity.  And
indeed we {\em cannot} say anything, because the rate of decay can be
arbitrarily slow.  What really is the content of the theorem then?

	The first basic observation is that the Fourier transform is
uniformly bounded under these conditions.  Namely,
\begin{equation}
\label{sup bound}
	\sup_{\xi \in {\bf R}} |\widehat{f}(\xi)| 
			\le \int_{\bf R} |f(x)| \, dx.
\end{equation}
This fact is very robust in terms of making finite approximations.

	The second observation is that for ``nice'' functions one can
get decay that is easily measured.  There is a general principle to
the effect that smoothness of a function is reflected in decay at
infinity of the Fourier transform, and there are many specific results
of this nature (often proved by integrating by parts).  Let us content
ourselves with a special case which is easy to calculate.  Let $a, b
\in {\bf R}$ be given, with $a < b$, and consider the function
\begin{equation}
	g_{a,b}(x) = \frac{1}{b-a} \, {\bf 1}_{[a,b]}(x).
\end{equation}
Here ${\bf 1}_A (x)$ denotes the characteristic (or indicator)
function of the set $A$, which means the function which takes the
value $1$ on the set and $0$ off the set.  The normalizing factor of
$\frac{1}{b-a}$ was included so that 
\begin{equation}
	\int_{\bf R} |g_{a,b}(x)| \, dx = 1.
\end{equation}

	The Fourier transform of this function is easy to compute.
It is given by
\begin{equation}
	\widehat{g}_{a,b}(\xi)  = \frac{1}{b-a} \, \frac{1}{- 2 \pi i \, \xi} 
		\, (e^{- 2 \pi i \, b \, \xi} - e^{- 2 \pi i \, a \, \xi}).
\end{equation}
Notice that this is compatible with (\ref{sup bound}) even if it may
not appear to be at first.  The potential for singularities provided
by the denominators is avoided by the difference in the exponentials.

	For all $a, b, \xi$ we have that 
\begin{equation}
\label{bound for the special case}
	|\widehat{g}_{a,b}(\xi)| \le \frac{2}{b-a} \, \frac{1}{2 \pi |\xi|}.
\end{equation}
That is, we have bounded the difference the exponentials by $2$.
For ``most'' $\xi$'s the difference in the exponentials will not be small.
It will be small only when $(b-a) \, \xi$ is close to an integer.

	This inequality makes concrete the second observation.  If $f$
is actually a step function -- a finite union of characteristic
functions of intervals -- then the Fourier transform decays at
infinity like a multiple of $\frac{1}{|\xi|}$.

	To prove the theorem one combines these two pieces of
information using a density argument, as follows.  Suppose that we are
given an integrable function $f(x)$ on ${\bf R}$, and let $\epsilon >
0$ be given.  A standard fact from real analysis is that there is a
step function $g$ on ${\bf R}$ such that
\begin{equation}
	\int_{\bf R} |f(x) - g(x)| \, dx < \epsilon.
\end{equation}
Using (\ref{sup bound}) we get that
\begin{equation}
	\sup_{\xi \in {\bf R}} |\widehat{f}(\xi) - \widehat{g}(\xi)| 
				< \epsilon.
\end{equation}
On the other hand our earlier calculations show that
\begin{equation}
	\lim_{|\xi| \to \infty} \widehat{g}(\xi) = 0.
\end{equation}
This implies that
\begin{equation}
	\limsup_{|\xi| \to \infty} |\widehat{f}(\xi)| \le \epsilon.
\end{equation}
Since $\epsilon > 0$ is arbitrary we conclude that the limit actually
vanishes.

	We cannot say anything about the rate at which the Fourier
transform decays at infinity.  There are two reasons for this.
Although we can approximate an integrable function $f$ by a step function
$g$ as above, we cannot say much about the nature of the approximation,
i.e., how will the complexity of $g$ grow as $\epsilon$ shrinks to $0$?
Without additional assumptions on $f$ the answer is that it can grow
as fast as it wants.  

	The second reason is that in our example above we had the
factor of $\frac{1}{b-a}$ in (\ref{bound for the special case}).  This
can be made as large as we want by taking $b$ close to $a$, without
affecting the $L^1$ norm of $g_{a,b}$.  We saw also that for most
$\xi$ we cannot really do much better than the estimate in (\ref{bound
for the special case}).  Of course we can build much nastier
integrable functions by combining infinitely many functions like
$g_{a,b}$ with $b-a$ as small as we want.  We have only to be careful
to include a summable series of coefficients to get an integrable
function in the end.

	The bottom line is that the rate of decay of the Fourier
transform of a general integrable function is unknowable because it
depends on the structure of a general integrable function, which can
be practically arbitrary.  In effect there are infinite processes in
both the statement and the conclusion of the theorem.

	We also saw another important point: that a crucial ingredient
for the result about the limit was the uniform bound (\ref{sup bound}),
which does have a more ``absolute'' meaning.  In particular it works
perfectly well in finite approximations.

\section{Lebesgue's theorem}
\label{Lebesgue's theorem}
\index{Lebesgue's theorem}

	Suppose that $f$ is an integrable function on ${\bf R}^n$.
Then
\begin{equation}
\label{lebesgue's theorem}
	\lim_{r \to 0} \frac{1}{|B(x,r)|} \, \int_{B(x,r)} f(y) \, dy
		= f(x)
\end{equation}
for almost all $x \in {\bf R}^n$.  Here $|B(x,r)|$ denotes the volume
of the ball $B(x,r)$.  Thus we are taking the average of the values of
$f$ on the ball of radius $r$ around $x$, then shrinking the radius to
$0$, and we are looking for the limit to be the original value at that
point.  This would occur at all points if $f$ were continuous, but for
functions which are merely integrable the limit need not even exist at
a given point.

	What can we really say about this limiting process?  What does
it mean concretely?  Is there any ``finite'' information contained
in this result?

	As in the previous section one cannot say anything about the
rate at which the limit converges.  One could not say anything even
for continuous functions.  This ambiguity is built into the statement,
because the continuity or integrability of a function does not give
terribly precise information about its structure.

	How can one prove such a theorem?  The basic structure of the
argument is much the same as before.  A key point is the density of
continuous functions.  Given an integrable function $f$ on ${\bf R}^n$
and an $\epsilon > 0$ there is a {\em continuous} function $g$
on ${\bf R}^n$ such that 
\begin{equation}
	\int_{{\bf R}^n} |f(x) - g(x)| \, dx < \epsilon.
\end{equation}
We want to use this to show that the existence of limits for
continuous functions implies the corresponding result for integrable
functions.

	To do this we need to have a uniform bound, analogous to
(\ref{sup bound}) in the previous situation.  In this case the
geometry is more tricky.  Consider the {\em maximal function}
\begin{equation}
	f^*(x) = \sup_{r > 0} \frac{1}{|B(x,r)|} \, \int_{B(x,r)} |f(y)| \, dy.
\end{equation}
In general when one is interested in the existence of a limit it is
helpful to first obtain uniform bounds.  This case is more complicated
because we need to take the supremum over $r$ for each point $x$, but
we need to treat the points $x$ more individually (and not just take
the supremum).

	Suppose for instance that we knew that there was a constant
$C > 0$ so that
\begin{equation}
\label{horrible lie}
	\int_{{\bf R}^n} f^*(x) \, dx \le C \, \int_{{\bf R}^n} f(x) \, dx 
\end{equation}
for all integrable functions $f$ on ${\bf R}^n$.  Then we would be
able to show that (\ref{lebesgue's theorem}) holds almost everywhere
for integrable functions using the fact that it is true for continuous
functions and an approximation argument, in much the same manner as
in the previous section.

	Unfortunately (\ref{horrible lie}) is wrong.  In fact, as soon
as $f$ is not $0$ almost everywhere, one has that
\begin{equation}
	f^*(x) \ge \frac{a}{1+|x|}
\end{equation}
for some constant $a > 0$ (which depends on $f$).  This is not hard to
verify, and it implies that
\begin{equation}
	\int_{{\bf R}^n} f^*(x) \, dx = \infty.
\end{equation}
The analogue of (\ref{horrible lie}) for $L^p$ spaces is true for all
$p > 1$, though.  There is a substitute for $p = 1$, which says that
$f^*$ is not too big most of the time.  The precise statement is that
there is a constant $C > 0$ such that
\begin{equation}
\label{weak type}
	\lambda \cdot |\{x \in {\bf R}^n : f^*(x) > \lambda \}| 
			\le C \, \int_{{\bf R}^n} f(x) \, dx 
\end{equation}
for all integrable functions $f$ on ${\bf R}^n$.  Note that (\ref{weak type})
would follow automatically if (\ref{horrible lie}) were true.

	It turns out that this weaker version is sufficient for
proving that (\ref{lebesgue's theorem}) holds almost everywhere for
integrable functions.  See \cite{Ga, St1, SW, Seapp}.

	There are several nice points now.  The first is that an
estimate like (\ref{weak type}) makes perfectly good sense in discrete
approximations.  In fact (\ref{weak type}) is equivalent to having
uniform bounds for its analogues in finite approximations.  The
theorem about the existence of limits does not directly make sense in
a finite context, but these bounds do make sense, and they are a large
part of the story.

	Actually there are ways to bring the existence of the limit
back to the finite world.  Although one cannot control the rate at
which the limit converges, one can count the number of oscillations.
That is, one can fix an $\epsilon > 0$ and ask how many times the
average
\begin{equation}
	\frac{1}{|B(x,r)|} \, \int_{B(x,r)} f(y) \, dy
\end{equation}
oscillates by more than $\epsilon$ as $r$ shrinks from $1$ down to
$0$, say.  If the limit as $r \to 0$ exists then the number of these
oscillations must be finite.  One can try to derive bounds for this
number, as a function of $x$.  Bounds on averages of these numbers, or
on the measure of the set of $x$'s where they are large, as in
(\ref{weak type}).  Such results can be derived using the methods of
Carleson's Corona construction, and results like this are explained
(in some versions) in \cite{Ga}.  Again these bounds are equivalent to
the existence of uniform bounds in finite contexts.

	Thus in fact there are ``finite'' statements which correspond
to existence almost everywhere of limits, even if there are problems
with the most obvious ways to do this.

	Another interesting aspect is that in many situations a
maximal estimate like (\ref{weak type}) is necessary to have existence
of limits almost everywhere.  That is, there are actual theorems to
the effect that existence of limits imply suitable bounds.  See
\cite{St, N, Mau}.  This fits nicely with the question of what really
exists at the level of ``finite'' mathematics.  It is typically easier
to make sense of the idea of a bound in finite contexts than the
existence of the limit, and although the two appear to be quite different,
in fact they are closely connected.

	There is a simpler and older result of this type in the
context of linear mappings between Banach spaces, called the uniform
boundedness principle.  \index{uniform boundedness principle}
(See \cite{Ru2}.)  Suppose that you have a
sequence of bounded linear mappings $T_j$ from some Banach space $(X,
\|\cdot\|_X)$ into another one $(Y,
\|\cdot\|_Y)$.  Suppose also that
\begin{equation}
\label{pointwise}
	\sup_j \|T_j(x)\|_Y < \infty 	\qquad\hbox{for all } x \in X.
\end{equation}
Then the $T_j$'s have uniformly bounded operator norms, which means
that there is a constant $M < \infty$ such that
\begin{equation}
\label{uniform bound}
	\sup_j \|T_j(x)\|_Y \le M \, \|x\|_X 	\qquad\hbox{for all } x \in X.
\end{equation}
Thus pointwise bounds imply uniform bounds, which is slightly
remarkable.  The point here is that $X$ is allowed to be
infinite-dimensional.  The theorem would be pretty trivial otherwise.

	Note that (\ref{pointwise}) holds automatically when
\begin{equation}
\label{pointwise limit}
	\lim_{j \to \infty} T_j(x) \hbox{ exists for all } x \in X.
\end{equation}
In practice it is not so easy to prove (\ref{pointwise limit}) 
without first establishing a uniform bound.

	The uniform boundedness principle provides another example of
a theorem which in effect seems to relate two kinds of infinite
processes in a nontrivial way while not appearing to have a natural
finite version.

\section{The ergodic theorem}

	There is a variation on the themes of the preceding section which
are closer to physical reality.  Let $(X, {\cal M}, \mu)$ be a finite
measure space, and let $T : X \to X$ be a measurable mapping which is
measure-preserving, so that
\begin{equation}
	\mu(T(A)) = \mu(A) \qquad\hbox{for all measurable sets } A \subseteq X.
\end{equation}
For example, think of $X$ as being the unit circle equipped with ordinary
Lebesgue measure, and let $T$ be rotation by the angle $\theta$.

	Fix a point $x \in X$, and consider the sequence
$\{T^j(x)\}_{j=1}^\infty$ of images of $x$ under $T$.  This is called
the {\em orbit} of $x$ under $T$.  How should we expect it to be
distributed in $X$?

	This question arises in classical mechanics.  Think of a
particle in a physical system moving under the influence of some
conservative force field.  For this one should think in terms of
continuous time instead of discrete time, as we have here, but this is
not a serious issue.  For instance, one can fix a small number $\delta
> 0$, and take $T$ to be the mapping that tells how a particle
described by $x$ moves after flowing in the system by an amount of
time equal to $\delta$.

	In the context of mechanics the assumption that $T$ be
measure-preserving is very natural, for reasons that we shall not
discuss.  (See \cite{Ar}.)  For this we should say that the phase
space $X$ tracks not only the position of the particle but also its
momentum.

	One might hope that the orbit of a point is uniformly
distributed throughout the space.  This is not true in general, and
one can make easy counterexamples.  It is true however for {\em almost
all} orbits when $T$ is an {\em ergodic} transformation.  This
assumption says that if $A
\subseteq X$ is measurable, and if $T(A) = A$, then either $\mu(A) =
0$ or $\mu(A) = \mu(X)$.  A subset $A$ such that $T(A) = A$ is said to
be {\em invariant}.  The existence of nontrivial invariant subsets
would prevent the orbits from being uniformly distributed, and so this
assumption is at least necessary.

	To understand what this means it is helpful to come back to
our example of a rotation on the unit circle by an angle $\theta$.  If
$\theta$ is a rational multiple of $\pi$ then the corresponding
rotation is not ergodic, and this is easy to check.  In this case the
mapping is even periodic, so that $T^n$ is the identity for some
positive integer $n$, and all orbits are finite.  It turns out that
the rotation is an ergodic transformation when the angle is an
irrational multiple of $\pi$.  See \cite{Sn} this and other more
difficult examples.

	The precise formulation of the {\em ergodic theorem}
\index{ergodic theorem} states that if $T : X \to X$ is
measure-preserving and ergodic, and if $f$ is an integrable function
on $X$, then
\begin{equation}
\label{ergodic theorem}
	\lim_{N \to \infty} \frac{1}{N} \sum_{j=1}^N f(T^j(x))
			= \int_X f \, d\mu
		\qquad\hbox{for $\mu$-almost all } x \in X.
\end{equation}
In other words, if we average $f$ over $X$, that is the same as
taking a random point $x \in X$ and averaging $f$ over the orbit
of $x$ under $T$.  Sometimes one says ``time averages equal space
averages''.  

	For example, we might take $f$ to be the characteristic
function of a measurable set $A$.  In this case the left side of
(\ref{ergodic theorem}) measures the amount of time that the orbit
of $x$ spends in $A$, while the right side reduces to the measure
of $A$.  Note that the existence of the limit on the left side
is part of the conclusion.

	This theorem is reminiscent of the existence of limits almost
everywhere in the preceding section, except that it goes in the
opposite direction, making averages in the large instead of over small
balls.  This turns out not to be a major issue.  For instance one of
the main ingredients in the proof is a estimate for the maximal
function
\begin{equation}
	\sup_N \frac{1}{N} \sum_{j=1}^N |f(T^j(x))|
\end{equation}
which is analogous to the one given before in (\ref{weak type}).

	See \cite{B} for a proof of the ergodic theorem and further
information.  See \cite{CW} for techniques for relating this more
closely to the matter of the previous section, through the method of
``transference''.

	If the mapping is measure-preserving but not ergodic one still
has bounds for the maximal function and the existence almost
everywhere of the time averages.  One simply cannot identify the time
averages with the space averages.

	As usual the existence of a limit does not have ``finite''
meaning but much of the mathematics around the theorem does have
finite meaning.  In this case however the physical interpretation
makes the question of finite meaning more compelling and perhaps more
intuitive or attractive.  For this one should look at more specific
situations and their dynamics, and this leads to a large area of
current mathematical research.

\section{Differentiability almost everywhere}
\label{Differentiability almost everywhere}

	Let $f : {\bf R}^n \to {\bf R}$.  We say that $f$ is differentiable
at the point $p$ if there is a linear mapping $A : {\bf R}^n \to {\bf R}$
(the differential of $f$ at $p$) such that
\begin{equation}
	f(x) = f(p) + A(x-p) + \epsilon(x),
\end{equation}
where 
\begin{equation}
	\lim_{x \to p} \frac{\epsilon(x)}{|x-p|} = 0.
\end{equation}

	There is a famous old theorem to the effect that Lipschitz
functions are differentiable almost everywhere.  (See \cite{Fe, St1,
Seapp}.)  As in Section \ref{Continuity}, a function $f : {\bf R}^n \to
{\bf R}$ is said to be {\em Lipschitz} if there is a constant $C > 0$
so that
\begin{equation}
\label{Lipschitz on R^n}
	|f(x) - f(y)| \le C \, |x-y| \qquad\hbox{for all } x, y \in {\bf R}^n.
\end{equation}

	This theorem is quite remarkable.  Although it can be seen as
a variant of the one discussed in Section \ref{Lebesgue's theorem}, it
is noticeably different.  In Section \ref{Lebesgue's theorem} we
started with an integrable function $f$ and we wanted to show that the
averages of the values of $f$ on small balls centered at a random
point $x$ converge to the value of $f$ at $x$.  In particular we had a
candidate for the limit in our hands already.  To prove that Lipschitz
functions are differentiable almost everywhere, part of the problem is
to find the candidate for the limit.

	Another point about this theorem is that Lipschitz functions
are arguably more ``elementary'' than integrable functions.  Just the
definition of Lebesgue integrals, integrability, and implicitly
measurability entail infinite processes of some sophistication.  In
the Lipschitz condition we have only to take a supremum over all pairs
of points and ask that the result be finite.  This is also an infinite
process, but it is more straightforward.

	The existence of a derivative of a Lipschitz function at a
point provides a kind of ``finiteness'' to its behavior.  It says that
the function is approximately affine at the given point, and affine
functions span a finite-dimensional space.  The asymptotic behavior of
Lipschitz functions on general metric spaces need not have ``finite
dimensional'' limitations on their asymptotic behavior of this type.

	Although one has the existence almost everywhere of the
derivative of a Lipschitz function, one cannot say anything about the
rate at which the convergence takes place.  How then might look for
natural versions of the theorem in finite contexts?  

	Normally we might ask for a uniform bound of the objects whose
limit we wish to take, but in this case we have simply assumed it in
(\ref{Lipschitz on R^n}).  (There is a different version of the story
in which uniform bounds on the difference quatients are not so immediate
and are more crucial to the existence of the limit, but we shall not pursue
that here.)

	Instead we can look for good bounds on the oscillation of the
function.  For instance, we can take a ball $B(x,t)$ and ask how well
$f$ is approximated by an affine function on this ball, as in
\begin{equation}
	\alpha(x,t) = \inf_A \, \sup_{y \in B(x,t)} \frac{|f(y) - A(y)|}{t},
\end{equation}
where the infimum is taken over all affine functions $A$.  For almost all
$x$ we have that
\begin{equation}
	\lim_{t \to 0} \alpha(x,t) = 0,
\end{equation}
because of the differentiability almost everywhere.  One can look for
quantitative bounds on the $\alpha(x,t)$'s which would make sense in
discrete approximations.  Such bounds exist, as in
\cite{Do}.  These bounds involve averages over $x$ and $t$, 
and not uniform rates of convergence.  

	In these bounds we would permit ourselves to approximate $f$
on balls by affine functions that are turning around as we move the
balls, but in fact one can also control the manner in which these
affine approximations are turning around.  This is similar to the
story of Carleson's Corona construction \cite{Ga}, as mentioned before
in Section \ref{Lebesgue's theorem}.

	These points about uniform bounds are discussed further in
\cite{DS4, Secon, Seapp}.

	The bottom line is that there are actually estimates which
work for discrete approximations and which capture rigidity properties
of Lipschitz mappings.

\chapter{Using infinite processes}

	Many times in mathematics it is very convenient to use
infinite processes, ones which are apparently much more complicated
than the problems to which they are applied.  We give some examples
here in this chapter, with emphasis on geometry and topology.  We are
particularly interested in situations where the statement is already
pretty concrete even if the proof is not.  The matter is not as much
fun when the statements themselves already rely on abstract processes.

	These examples should be compared to the discussion of lengths
of proofs in Chapter \ref{About symmetry in language}.

\section{Mostow rigidity}

	We consider a rigidity phenomenon which is closer to
differential geometry and Lie groups.  For this discussion it will be
convenient to have the idea of {\em hyperbolic geometry}.  The easiest
way to describe this is to use a basic model, the upper half-space
${\bf R}^n_+$, which means the set of points whose last coordinate is
positive:
\begin{equation}
	{\bf R}^n_+ = \{x \in {\bf R}^n : x = (x_1, \ldots, x_n), \, x_n > 0\}.
\end{equation}
We can define a {\em Riemannian metric} on this set using
\begin{equation}
	\frac{dx_1^2 + \cdots + dx_n^2}{x_n^2}.
\end{equation}
One can think of this simply as a recipe for computing infinitesimal
lengths of curves.  The lengths are distorted by a factor of $x_n^{-1}$
as we move towards the boundary hyperplane $x_n = 0$.  ``Distorted''
means relative to ordinary Euclidean geometry.

	To give a feeling for the behavior of this geometry, consider
the transformation
\begin{equation}
	x \mapsto a \, x, 
\end{equation}
where $a$ is a positive real number.  This transformation has the
simple effect on Euclidean geometry of changing distances by exactly a
factor of $a$.  As for our hyperbolic metric, notice that
transformations of this form map ${\bf R}^n_+$ onto itself, and that
{\em they do not change hyperbolic distances}.  One can use this
together with the obvious invariance of the metric under translations
$x \mapsto x + b$, $b \in {\bf R}^n$, $b_n = 0$, to show that {\em all
points in ${\bf R}^n_+$ look the same for the hyperbolic metric}.
That is, given any pair of points in ${\bf R}^n_+$, one can find a
mapping from ${\bf R}^n_+$ onto itself which preserves the hyperbolic
geometry.

	We shall restrict ourselves to $n \ge 2$ for this discussion.
Although the definition makes sense for $n=1$, it does not give back
anything new, just another copy of the Euclidean line.  The mapping
$x \mapsto e^x$ provides the correspondence between the two metrics.

	We can define hyperbolic geometry as a general concept using
this model.  Let $M$ be a smooth manifold of dimension $n$.  If we
have chosen a Riemannian metric on $M$, then we call this metric {\em
hyperbolic} if around each point in $M$ there is a coordinate chart
which identifies a neighborhood of the point in $M$ with a
neighborhood of a point in ${\bf R}^n_+$ in such a way that the metric
on $M$ corresponds exactly to our metric on ${\bf R}^n_+$.  There
are various ways to formulate this, in terms of ``constant negative
curvature'', for instance, but they amount to this in the end.

	An amazing theorem of Mostow \cite{M} implies that if $M_1$
and $M_2$ are two compact hyperbolic manifolds of the same dimension
$n$ whose fundamental groups are isomorphic, and if $n > 2$, then
$M_1$ and $M_2$ are {\em isometrically equivalent}.  This means that
there is a mapping from $M_1$ onto $M_2$ which preserves the metric.
This is extremely wrong in dimension $2$, in which surfaces can admit
{\em continuous} families of inequivalent hyperbolic geometries.

	What does this theorem mean in terms of ``finite'' mathematics?

	To understand this question better it is helpful to consider a
reformulation.  If $M$ is a compact hyperbolic manifold of dimension
$n$, then it can be realized as a quotient of our model space ${\bf
R}^n_+$.  That is, there is a discrete group $\Gamma$ of isometries of
${\bf R}^n_+$ such that $M$ is isometrically equivalent to ${\bf
R}^n_+ / \Gamma$.  The {\em isometries} of ${\bf R}^n_+$ are the
mappings from ${\bf R}^n_+$ onto itself which preserve the hyperbolic
metric.  This is a group, and in fact it can be described as a
classical Lie group of matrices.  It can also be realized as the group
generated by translations and dilations mentioned above together with
reflections about spheres whose centers lie on the hyperplane $x_n =
0$.

	This representation of $M$ as ${\bf R}^n_+ / \Gamma$ is not
terribly mysterious.  Another way to say it is that the universal
covering of $M$ is isometrically equivalent to ${\bf R}^n_+$.  By
assumption $M$ is locally equivalent to ${\bf R}^n_+$, and so the
point is to knit together many local equivalences into a global one.
By passing to the universal covering we avoid ambiguities that could
come from having a sequence of local equivalences loop around to a
place where they have been before, but where a different choice of
local equivalence was used.

	At any rate this representation is well known.  The group
$\Gamma$ is isomorphic to the fundamental group of $M$ as a discrete
group.  Thus in the context of Mostow's theorem we can represent $M_1$
and $M_2$ as ${\bf R}^n_+ / \Gamma_1$ and ${\bf R}^n_+ / \Gamma_2$,
where $\Gamma_1$ and $\Gamma_2$ are groups of isometries on ${\bf
R}^n_+$ which are isomorphic as groups to the fundamental groups of
$M_1$ and $M_2$.  These fundamental groups are isomorphic to each
other by assumption, and so $\Gamma_1$ and $\Gamma_2$ are isomorphic
to each other.

	Isomorphic as abstract groups, that is.  That does not
necessarily mean that they sit inside of the group of all hyperbolic
isometries on ${\bf R}^n_+$ in the same way.  Mostow's theorem implies
that they do.  The isometry between $M_1$ and $M_2$ comes from an
isometry on ${\bf R}^n_+$, and this implies that $\Gamma_1$ and
$\Gamma_2$ are actually conjugate to each other within the group of
isometries on ${\bf R}^n_+$.  This can fail when $n = 2$ or for
discrete groups of isometries which do not correspond to compact
manifolds, although the compactness assumption can be weakened.

	This rigidity phenomenon is actually part of a vast area in
mathematics, involving other kinds of groups and geometries.  See
\cite{GP} for a survey.

	There is an amusing aspect of Mostow's original proof.  It is
not so difficult to make the groups $\Gamma_1$ and $\Gamma_2$
correspond under a mapping from ${\bf R}^n_+$ onto itself, the problem
is that this mapping need not be an isometry.  This mapping can be
taken to be nice and smooth on ${\bf R}^n_+$, but it is more
interesting to look at its behavior at the boundary hyperplane.  At
the boundary the mapping need not be smooth a priori, but it does
satisfy a geometric condition called ``quasisymmetry'', which means
that it distorts relative distances in a bounded way.  It turns out
that this condition implies the existence of differentials at almost
all points on the boundary, differentials that are nondegenerate.
This fact can be seen as a relative of the differentiability almost
everywhere of Lipschitz mappings.  (It is important here that our
initial dimension be larger than $2$, so that the dimension of the
boundary hyperplane be at least $2$.)  The existence of a well-behaved
differential at a point on the boundary provides a bit of asymptotic
flatness which can then be ``blown up'' to the whole space using the
invariance of the problem under our groups $\Gamma_1$ and $\Gamma_2$.
More precisely one can show that our mapping on ${\bf R}^n_+$ which
which we know cooperates with the groups $\Gamma_1$ and $\Gamma_2$ but
is not necessarily an isometry can be deformed into a mapping which
does both.  This uses also the compactness of $M_1$ and $M_2$ to
ensure that the groups $\Gamma_1$ and $\Gamma_2$ are large enough that
they both ``see'' all of ${\bf R}^n_+$, and in particular they can
detect the point of asymptotic flatness on the boundary.

	This outline is vague but it brings up an interesting point.
In the middle of the proof one uses a theorem about differentiability
almost everywhere to get the existence of a certain kind of limit even
though the final statement does not really call for such a limit.  One
might say that it is a more infinite process than the statement seems
to need.  This is particularly true when one thinks about starting
with finite (smooth) hyperbolic manifolds, which are more like
``finite'' complexity.  The mapping that we create at the boundary is
more like ``infinite'' complexity, although it has a huge amount of
symmetry coming from the way that it interacts with the groups
$\Gamma_1$ and $\Gamma_2$.  In the step where we obtain
differentiability at some point we forget about the extra symmetry for
a moment and simply treat the mapping as a general one.

	Thus {\em this} method, which was the original one, does seem
to leave the ``category'' of complexity of the original statement in a
rather clear and interesting way.

	See \cite{GP} for other approaches to rigidity.

\section{Groups of polynomial growth}

	Let $G$ be a finitely generated group with generating set $S$,
say.  Let $V(n)$ denote the number of elements of $G$ which can be expressed
as a word of length $\le n$ in the generating set $S$.  We say that
$G$ has {\em polynomial growth} if there are constants $C, d > 0$
such that
\begin{equation}
	V(n) \le C \, n^d \qquad\hbox{ for all } n \le 1.
\end{equation}
This property does not depend on the particular choice of (finite)
generating set $S$.

	Gromov \cite{Gr} proved that a finitely generated group of
polynomial growth contains a nilpotent group of finite index.  The
proof was more transcendental than one might expect the statement to
require.  Roughly speaking one builds certain auxiliary ``continuous''
spaces on which the group acts by isometries.  These auxiliary
continuous spaces are obtained naively from the group as limits of
discrete spaces.  The isometry groups of these auxiliary continuous
spaces are shown to be Lie groups (with only finitely many components)
through general results of Montgomery and Zippin.  This method permits
one to obtain a lot of homomorphisms into Lie groups, which is an
important part of the argument.

	In this case the statement of the theorem is already pretty
concrete.  The assumption of polynomial growth entails an infinite
process, but it is not so bad.  The method of the proof relies on
processes which are much more transcendental.

\section{A product expansion for $\sin z$}

	Consider the formula
\begin{equation}
\label{identity}
	\sin \pi z = 
		\pi z \prod_{n=1}^\infty \biggl(1 - \frac{z^2}{n^2}\biggr).
\end{equation}
The right-hand side makes sense because of the summability of
$\sum \frac{1}{n^2}$.

	This is a well-known identity from complex analysis, as on
p195 of \cite{A}.  The proof is fairly easy to understand.  We think
of both sides of the equation as being functions of the complex
variable $z$ on the whole complex plane.  It is a standard fact that
the convergence of the product is sufficiently tame to give a {\em
holomorphic} function on the plane.  The product was chosen so that
both sides of the equation vanish at exactly the same points, and so
that the product is the ``simplest'' function with these properties.
This permits us to take the quotient of the two functions to get a
holomorphic function on the whole complex plane which never vanishes.
The point then is to show that this quotient is actually constant.  If
it is constant then it must be equal to $1$, because the quotient
takes the value $1$ at the origin.

	For the proof that the quotient is constant we refer to \cite{A}.
Instead of the proof given on p195 I prefer to make an argument based
on Theorem 8 on p207 of \cite{A}.  The idea is that one has some control
on the behavior of our functions at infinity, in such a way that we can
conclude that the logarithm of the quotient is linear.  The quotient
is an even function by inspection, which permits us to conclude that
the quotient is actually constant.

	Thus at bottom we derive our identity from the standard
miracle of complex analysis, which is that entire holomorphic
functions of moderate growth are very special.  These are variations
on the theme of Liouville's theorem, to the effect that a bounded
entire holomorphic function must be constant.  The proof given on p195
of \cite{A} uses a similar argument.

	On the other hand we can see (\ref{identity}) as a sequence of
identities about infinite series.  We know the Taylor expansion of
$\sin \pi z$,
\begin{equation}
    \sin \pi z = \sum_{n=0}^\infty (-1)^n \, \frac{(\pi z)^{2n+1}}{(2n+1)!}.
\end{equation}
One can also compute the Taylor coefficients of the right side of
(\ref{identity}), simply by multiplying out the product and collecting
terms to get the Taylor series expansion about the origin.

	Of course (\ref{identity}) is equivalent to the equality of
the two Taylor expansions.  The individual equalities between Taylor
coefficients are not quite ``finite'' statements in mathematics, since
they involve $\pi$ and infinite sums, but they are fairly concrete.
They have natural finite approximations, and all of the relevant
series converge absolutely.

	Thus again we have a situation where the mathematical content
of the statement is arguably much less ``infinite'' or abstract than
the proof through complex analysis.  This case is much easier to
understand than the story of groups of polynomial growth.

\section{Transversality}
\index{transversality}

	The following is a version of ``Sard's theorem''.

\begintheorem
Suppose that $U$ is an open subset of some ${\bf R}^m$ and that
$f : U \to {\bf R}^n$ is smooth.  Set
\begin{equation}
	{\cal C} = \{x \in U : {\rm rank} \, df_x < n \}.
\end{equation}
Then the image of ${\cal C}$ under $f$ has Lebesgue measure zero
in ${\bf R}^n$.
\end{theorem}

	Here $df_x$ denotes the {\em differential} of $f$ at $x$, in
the sense of advanced calculus.

	See \cite{Fe, Mi2} for proofs and more information.  Note that
the theorem is more subtle when $n < m$ than when $m = n$.

	A particularly useful consequence of this theorem is that
${\bf R}^n \backslash f({\cal C})$ is everywhere dense, and this was
proved earlier by Arthur Brown.  See \cite{Mi2, MiS} for further
information about the history of this result.

	In topology one sometimes uses these results for the sake of
{\em transversality}.  Here is an example.

\begintheorem
\label{f^{-1}(p)}
Suppose that $M$ and $N$ are smooth manifolds of dimension $m$ and $n$,
respectively, with $n \le m$, and suppose that $f : M \to N$
is a smooth mapping.  Then $f^{-1}(p)$ is a smooth embedded
submanifold of $M$ for ``almost all'' $p \in N$.
\end{theorem}

	As usual, ``almost all'', meaning all but a set of (Lebesgue)
measure zero, makes sense for any smooth manifold, without extra
structure.  In order to speak about the actual {\em measure} of a set
-- and not simply whether or not the measure is zero -- requires
additional structure, such as a volume form, or a Riemannian metric.

	This theorem follows easily from the preceding one and the
implicit function theorem.  For a discussion of the latter in the
context of smooth manifolds see \cite{Wa}.

	This is a kind of ``transversality'' theorem.  It is a basic
tool in topology for demonstrating the existence of manifolds or
submanifolds with prescribed properties while not {\em constructing}
them exactly, and certainly not through an explicit finite recipe.
For instance, it plays a role in the assertion that if $M$ is a
compact manifold of dimension $m$ and $0 \le k \le m$ is an integer
which satisfies $m < 2
\, k -1$, then the  {\em rational} cohomology of $M$ in dimension $k$
is generated by the Poincar\'e duals of {\em smooth} submanifolds of
$M$ of dimension $m-k$ whose normal bundles are {\em trivial}.  See
the discussion on p232-233 in \cite{MiS}.  (It is important to exclude
torsion in this discussion.)

	The celebrated work of Thom \cite{Th} provides an analysis of
when a given homology class is represented by a smmoth submanifold, or
when a given manifold can be realized as the boundary of another.
These results rely heavily on transversality arguments, but of a
slightly different nature from the ones above.  Given a smooth mapping
$f : M \to N$ and a smooth submanifold $Y$ of $N$, one is interested
to know if $f$ is {\em transverse} to $Y$, or can be made so after
small perturbation, and such {\em transversality} would imply that
$f^{-1}(Y)$ is a smooth submanifold of $M$ of the correct dimension.
See \cite{MiS} for more information.

	In this way one might obtain smooth manifolds or classes of
smooth manifolds with certain properties, but in terms of {\em
complexity} this method has the unfortunate feature that it provides
little information of a quantitative nature.  In the context of
Theorem \ref{f^{-1}(p)}, for instance, one might be able to get bounds
on the {\em mass} of the submanifolds obtained through the {\em
co-area theorem}, as in \cite{Fe}, which provides bounds on the
``average'' mass of $f^{-1}(p)$ for $p \in N$, but in general one
might be hard pressed to do better.  Even for that one would need
bounds on the size of the differential of $f$, and they might not be
so easy to obtain.

	For a broader discussion of geometric bounds in the context of
algebraic topology see \cite{Gr0}.  Indeed it was Gromov who first
pointed out to me that there are basic questions of complexity to be
addressed concerning the cobordism theory of Thom, and also the role
of complexity and Lipschitz bounds for homotopy classes of mappings
between manifolds (especially in connection with classical results of
Serre).

	Note that these issues can be discussed in the ``finite''
contexts of polyhedra instead of manifolds.  For a combinatorial
version of Theorem \ref{f^{-1}(p)}, for instance, see p235 of
\cite{MiS}.  In this case transversality can be seen as finite
and combinatorial.

\section{Homeomorphisms}
\label{Homeomorphisms}

	Suppose that $P$ and $Q$ are finite polyhedra.  One can think
of them as being finite unions of simplices in some Euclidean space,
or in purely finite terms, characterized by the combinatorics of the
intersections of the sub-simplices.  There is a natural notion of {\em
combinatorial equivalence}, in which we ask that there be a
correspondence between the simplices inside $P$ and $Q$ which respects
the way that they meet in faces, or that this be true after we {\em
subdivide} $P$ and $Q$.  The effect of allowing subdivisions is to
treat a single line segment and the union of two segments which meet
exactly in an endpoint as being the same, for instance.

	One can think of combinatorial equivalence geometrically in
terms of the existence of a homeomorphism between $P$ and $Q$ which is
{\em piecewise-linear}.  For this we use the linear structure from the
ambient Euclidean space, and the fact that simplices are truly flat.

	There was a famous question (the ``Hauptvermutung'') in
topology about whether two finite polyhedra which are homeomorphic
are actually combinatorially equivalent.  This is a very nice
question, because general homeomorphisms allow infinite processes,
while combinatorial equivalence is finite.  As a practical matter a
positive answer would be useful because it would allow one to show
that combinatorial invariants of polyhedra lead to topological
invariants.  Classically one was interested in the topological
invariance of homology groups defined combinatorially, but this
was resolved in another way.

	A counterexample to this question was found by Milnor
\cite{Mi}.  A more striking example appeared later from work of Edwards
\cite{E}, who found exotic triangulations of spheres.  
See also \cite{C1, C2, Dv1}.  For this we take $Q$ to be a standard
polyhedral representation of the $5$-sphere, e.g., as the boundary of
a $6$-dimensional simplex.  The statement then is that there is
another $5$-dimensional finite polyhedron $P$ which is homeomorphic to
${\bf S}^5$ but not piecewise-linearly equivalent to $Q$.  This
polyhedron $P$ is obtained in a specific way, as the ``double
suspension of a homology $3$-sphere'', which forces any homeomorphism
between $P$ and $Q$ to undergo enormous stretching.  More precisely,
there is a particular polygonal curve $\Gamma$ in $P$ whose complement
in $P$ is not simply connected, and the image of this curve in $Q$
under any homeomorphism must have ``Hausdorff dimension at least
$3$'', as pointed out in \cite{SS}.  The latter condition implies in
particular that if $h : P \to Q$ is a homeomorphism, then for
sufficiently small $\epsilon > 0$ we cannot cover $h(\Gamma)$ by fewer
than about $\epsilon^{-3}$ balls of radius $\epsilon$.  By contrast,
for $\Gamma$ itself we can use about $\epsilon^{-1}$ such balls, since
it is polygonal (a finite union of straight line segments).

	This implies that a homeomorphism from $P$ onto $Q$ cannot
even satisfy a H\"older condition like (\ref{Holder}) in Section
\ref{Continuity} unless the exponent $\alpha$ is at least $1/3$.  It
is not known whether there are homeomorphisms in this case which are
H\"older continuous of any positive exponent $\alpha$.

	So what does this mean in terms of finite versus infinite
processes in mathematics?  What would people have thought earlier in
the century when they were debating about foundations and mathematical
logic if they knew about {\em this}?

	A slightly different type of situation occurs in dimension
$4$.  There are compact smooth manifolds of dimension $4$ which are
homeomorphic to each other but not diffeomorphic.  Because these
manifolds are smooth we can realize them as finite polyhedra in a way
that is compatible with their smooth structures.  These polyhedra are
not exotic like Edwards' example; locally they are equivalent to
standard (flat) $4$-space.  If the smooth manifolds were diffeomorphic
to each other then the corresponding polyhedra would be
combinatorially equivalent.  In this context they are not
combinatorially equivalent, and the homeomorphisms between them must
distort distances in a strong way.  The necessity of unbounded
distortion is not as easy to see in this situation as for Edwards'
example.  It is shown in \cite{DoS} that the homeomorphisms cannot be
``bilipschitz'' or even ``quasiconformal''.  The latter means that
even {\em relative} distances have to be strongly distorted.  (See
\cite{DoS} for details.)  It is not known whether there exist
homeomorphisms with bounded distortion of relative distances in the
context of Edwards' example and its extensions.

	The existence of homeomorphisms in these examples of
$4$-dimensional manifolds is provided by the work of Freedman
\cite{FQ}, whose constructions entail very complicated infinite
processes.

	This type of phenomenon does not occur in dimensions less than
$4$.  For dimensions $\le 3$ it is true that the existence of a
homeomorphism between two smooth manifolds implies the existence of a
diffeomorphism.  In dimensions larger than $4$ it can again happen
that compact smooth manifolds are homeomorphic but not diffeomorphic
or even piecewise-linearly equivalent.  (See \cite{KS, Si}.)  In
dimensions $\ge 5$ however there will always then exist a
homeomorphism which is ``bilipschitz'', and therefore does not distort
distances by more than a bounded amount.  This was proved by Sullivan
\cite{Su}.

	Thus for {\em manifolds} of dimensions $\ge 5$ there are still
infinite processes involved but they are apparently much milder than
the ones in dimension $4$.  Still their precise nature remains unclear
though.

	These examples from topology are very striking.  Their
geometric concreteness -- finite polyhedra! -- lead us to hope for
similar concreteness for the associated mappings, but this simply is
not true.

	What are we to make of this then?  Are these homeomorphisms so
bizarre that we should reject them as being too far from concreteness?
What would people have thought earlier in the century?

	Unfortunately we do not understand so well how these
homeomorphisms look.  A better picture of the infinite processes
involved, of natural finite approximations to them, or algorithmic
processes for calculting them, would make them more concrete.

	In the case of manifolds of high dimension there are fairly
concrete examples \cite{Si}.  For the strange polyhedral spheres of
Edwards and Cannon and the strange mappings and spaces in dimension
$4$ of Freedman the picture is much less clear.  It can be easy to
understand the spaces concretely, in finite terms perhaps, but the 
mappings between them are another matter.

	There is a lot of room for better theories of complexity in
geometry and topology.  A basic issue is to find natural phenomena of
``finiteness'' for the complexity of spaces and mappings between them.
A first point for this is that there are only countably many compact
topological manifolds up to homeomorphism type \cite{CK}.  This
suggests that there is a possibility for a more ``concrete''
understanding.

	There has been a lot of work done on finding geometric conditions
on manifolds which imply that there are only finitely many underlying
topological types.  See \cite{Ch, F2, GPW, P1, P2}, for instance.

	In \cite{F1} there is a result to the effect that a certain
kind of almost-homeomorphism is practically as good as a
homeomorphism, at least in high dimensions.  This type of result is
very attractive if one wants to find ``finite'' formulations of
otherwise infinite notions in topology.  (One should be careful about
the details here, about the quantifiers.  Depending on how one
formulates a notion of almost-homeomorphisms, how much one is willing
to assume, it may be more or less easy to show that they can be
deformed into homeomorphisms.  If one makes very strong assumptions
this can be achieved through a simple compactness argument.  The result in
\cite{F1} is much better than that.)

	The aforementioned finiteness and countability results involve
somewhat similar considerations (e.g, notions of spaces being almost
homeomorphic).

	Alternatively one can try to find better ways to manage
infinite degrees of complexity, in particular to measure them in
intelligent ways.  This is more the perspective of \cite{Secon,
Seapp}.  The idea of infinite complexity is quite basic to analysis,
even if it is not described explicitly as such.  This idea fits
particularly well with the discussion of Section
\ref{Differentiability almost everywhere}.  See \cite{DS4} for related
issues of quantitative geometry and infinite complexity (some of which
are also explained in the more expository \cite{Secon, Seapp}).

	The idea of complexity in geometry and topology is much
broader than the particular areas mentioned here.  See \cite{Gr0}, for
instance.

\chapter{Finite constructions and their limits}

	The idea of continuous processes is very powerful in
mathematics.  As in calculus, for instance.  On the other hand the
reliance of mathematics on continuous processes may lead one to
overlook many interesting phenomena.  If one thinks in terms of
continuous processes one is likely to see continuous processes.

	One can see the limitations of the continuous realm in all
sorts of combinatorial problems.  In this chapter we try to go closer
to the edge between them.

	This is a murky area.  We are not ready to formulate a clear
general principle or anything like that.  We shall often find it
necessary to abuse the small amount of general language that we have.
In one moment we may wish to consider the real numbers as a
``continuous'' object, certainly so in contrast to the integers or the
rationals, while in another moment we may wish to view them as
``finite'', only one dimension, the ``continuous'' versions to be
found in infinite-dimensional spaces.

	Still one can attempt some general language.  Often in
mathematics one is interested in objects which are somehow finite but
large, a finite number of particles, a finite number of parameters.
One is interested in behavior as the number goes to infinity.
Sometimes there are natural limiting objects to consider, sometimes
not.  Roughly speaking we shall try to aim towards situations where
one might wish that there were a continuous limit but there is not a
natural one, or where continuous models do not capture all that one
wants to analyze.

	It is healthy to keep in mind that for purposes of computation
or physics one is often interested in numbers which are large but
``far'' from infinite.  Numbers of modest logarithm, for instance.
The restriction to relatively ``concrete'' large numbers is sometimes
modelled well by the idea of asymptotic behavior at infinity, and
sometimes not.  Sometimes one cares about what happens at the level of
$10^{10}$ but not at the level of $10^{20}$, and this can effect the
relevant mathematical structure.  Let us begin with an example to
illustrate this point.

\section{The Hilbert transform}
\label{The Hilbert transform}
\index{Hilbert transform}

	Let $f$ be a function on the real line, and let us define a
linear operator $H$ on such functions by
\begin{equation}
	H f(x) = {\rm p. v.} \int_{-\infty}^{\infty} \frac{f(y)}{x-y} \, dy.
\end{equation}
The existence of this integral is problematic because of the
singularity, but imagine writing it as
\begin{equation}
\label{computation}
	\int_{|x-y| > 1} \, \frac{f(y)}{x-y} \, dy +
  		\lim_{\epsilon \to 0} \, \int_{\epsilon < |x-y| \le 1} 
					\, \frac{f(y) - f(x)}{x-y} \, dy.
\end{equation}
We can subtract off the $f(x)$ in the last term because
\begin{equation}
\label{cancellation}
	\int_{\epsilon < |x-y| \le 1} \frac{1}{x-y} \, dy = 0,
\end{equation}
as one can easily check, using the fact that $\frac{1}{x-y}$ is
antisymmetric about $x$.  This computation implies that the principal
value which defines $H f(x)$ exists whenever $f$ has compact support
and is differentiable at $x$, for instance.  One can do much better
than this, and there are deeper results to the effect that $H f(x)$
exists almost everywhere when $f \in L^p({\bf R})$, $1 \le p <
\infty$.  See
\cite{Ga, Ka, St1}.

	This operator is called the {\em Hilbert transform}.  It is very
basic in analysis and is connected to complex analysis, partial differential
equations, and operator theory.  It turns out that it actually defines
an isometry on $L^2({\bf R})$, modulo a constant multiplicative factor.
One can see this using the {\em Fourier transform}, for instance.
A deeper fact is that $H$ defines a bounded linear operator on
$L^p({\bf R})$ for each $1 < p < \infty$.  See \cite{Ga, Ka, St1}.

	Thus $H$ does not distort the size of a function too badly.
This is a little surprising if one simply stares at the integral
naively.  For instance, what if one tries to put in absolute value
signs?  What about
\begin{equation}
	\int_{-\infty}^{\infty} \frac{|f(y)|}{|x-y|} \, dy?
\end{equation}
This is a disaster.  It will be equal to $+ \infty$ at almost all
points $x$ at which $f(x)$ is nonzero (at all of them if $f$ is
continuous).  Because we put in the absolute values we cannot make a
trick like (\ref{computation}).  The analogue of (\ref{cancellation})
does not work now.

	Thus the ``cancellation'' in the kernel of the Hilbert
transform plays a crucial role.  Without the cancellation the operator
is hopeless, the integral simply diverges.  

	But is the bad operator really so bad?  To understand this it
is helpful to think about what happens when one makes discrete
approximations.  The fact that one has bounds for $H$ -- e.g., that it
defines a bounded operator on $L^p$ -- is essentially the same as
saying that discrete approximations to $H$ satisfy uniform bounds,
independently of the level of the discretization.  For the bad
operator, with the absolute values inside, the discretizations will
become unbounded in the limit.

	Let us be more precise.  Imagine replacing the real line with
the set $X_N$ of numbers between $-N$ and $N$ which are integer
multiples of $\frac{1}{N}$, where $N$ is a large integer.  It is easy
to define versions of $H$ and its ``bad'' counterpart with the
absolute values inside on the discrete model $X_N$.  A function is now
a function on $X_N$, and one replaces integrals and integral norms
with sums.  The analogue of $H$ on $X_N$ has $L^p$ bounds which do not
depend on $N$ when $1 < p < \infty$.  One can show that the analogous
operator with the absolute values inside has operator norm on $L^p$
like $O(\log N)$.

	For computationally realistic values of $N$ the ``bad''
operator is essentially bounded.  This was first pointed out to me by
R. Coifman, who got it from V. Rokhlin.  For me this was quite
startling, something to resist at first, but indeed the ``bad''
operator is not really so bad computationally.  It can be computed
quickly by a computer using wavelet-based methods, as in \cite{BCR}.

	In summary, although boundedness properties of singular
integral operators correspond exactly to uniform bounds for their
discrete approximations, some operators that are bad in the limit are
pretty good for realistic levels of discretization.

\section{Topology}

	It is very natural to ask about the {\em topological}
equivalence of two spaces.  To try to track the way that 
pieces are connected together if not their precise size and shape.

	In some cases we can try to think about this in finite terms.
For instance we can look at {\em finite polyhedra}, and we can look at
{\em combinatorial} or (which is the same) {\em piecewise-linear}
equivalence.  These are finite notions which are basically
topological.

	We also have the usual notion of topological equivalence,
defined in terms of the existence of a homeomorphism.  Topological
equivalence is implied by combinatorial equivalence, but the converse
can fail dramatically, as we discussed in Section
\ref{Homeomorphisms}.

	Thus in topological equivalence we have a concept which is a
kind of limiting version of combinatorial equivalence but which does
not reduce to the latter in the cases where it could, i.e., finite
polyhedra.  One can try to {\em make} finite versions of it for finite
polyhedra, through mappings which might behave like a homeomorphism at
all scales above some particular small scale.  This gets to be
complicated, especially if we are not specifying a rate of continuity
for our mappings in advance.  In the infinite world we can use the
notion of continuity to accommodate all rates of continuity at once,
as in Chapter \ref{Examples from the infinite world}, while in the
finite world this does not work so easily.

	If we are simply interested in {\em mappings}, rather than
{\em equivalences}, then the story is simpler.  We might view
continuous mappings between topological spaces as an ``infinite''
version of continuous piecewise-linear mappings between finite
polyhedra.  The ``completion'' to the infinite notion behaves fairly
well, in the sense that arbitrary continuous mappings between finite
polyhedra can be approximated by piecewise linear mappings.  This is
well-known and not difficult to prove, and its analogue for
homeomorphisms is far from true, as we have seen.
	
	Consider now the notion of {\em homotopy equivalence}.
\index{homotopy equivalence}  Suppose that we are given two finite
polyhedra $P$ and $Q$, and a pair of continuous maps $f, g : P \to Q$.
We can talk about the mappings being homotopic, so that there is a
continuous mapping $H : P \times [0,1] \to Q$ which interpolates
between them, $H(x,0) = f(x)$, $H(x,1) = g(x)$.  If we have continuous
mappings $f : P \to Q$ and $\phi : Q \to P$ such that $\phi \circ f :
P \to P$ and $f \circ \phi : Q \to Q$ are both {\em homotopic} to the
identity mappings on $P$ and $Q$, then we get a {\em homotopy
equivalence} between $P$ and $Q$.  This defines another topological
notion of equivalence between finite polyhedra, which makes sense as
well for topological spaces in general.  While this notion lives a
priori in the infinite world, since we are allowing continuous
mappings, we could also define a ``finite'' version, in which we
restrict ourselves to mappings which are piecewise linear.  This would
include the ``interpolating'' mappings implicit in the homotopies
(like $H$ above).

	These two versions of homotopy equivalence actually define the
same relation between finite polyhedra.  This is well known and not
hard to prove, using the fact that continuous mappings between finite
polyhedra can be approximated (uniformly) by piecewise linear
mappings, as mentioned above.  Thus we have a simpler correspondence
between the ``finite'' and ``infinite'' worlds for the notion of
homotopy equivalence than we have for topological equivalence.

	There are variations on this theme, in which we look at
mappings which are more than continuous but less than
piecewise-linear, but the matter is about the same.  One could use
only continuous mappings which are Lipschitz or H\"older continuous,
for instance.  In that case the notion of homotopy equivalence would
not apply to topological spaces in general but to metric spaces, for
instance.

	This brings up an interesting point.  Suppose that we begin
with the idea of homotopy equivalence between finite polyhedra in
which all the mappings are piecewise linear, and imagine that we ask
ourselves naively for natural ``infinite'' versions.  We might find
ourselves lead to homotopy equivalence between arbitrary topological
spaces, or to homotopy equivalence between metric spaces in which all
the mappings are Lipschitz or H\"older continuous, or to something
else besides.  Thus we can have different ``completions'' of the
original finite concepts.  Different but rather close, enjoying
natural compatibilities.

	If we were to start with combinatorial equivalence and look
for natural ``continuous'' versions, we might arrive at the notion of
homeomorphic equivalence between topological spaces, or something
stronger for metric spaces, involving Lipschitz of H\"older
continuity, for instance.  Again we get different ``completions'' to
the infinite world.  This time they are much more different, without
such simple compatibilities as before.

	Thus finite notions might or might not have well-behaved
continuous versions, and these versions might or might not be
``unique''.

	(Some additional technical remarks.  Let us say that two
metric spaces are {\em bilipschitz equivalent} \index{bilipschitz
equivalent} if there is a bijection between them which does not expand
or contract distances by more than a bounded factor.  This is the same
as asking that both the bijection and its inverse be Lipschitz
mappings.  Combinatorial equivalence between finite polyhedra implies
bilipschitz equivalence.  We already know that topological equivalence
is strictly weaker than combinatorial equivalence for finite
polyhedra, and bilipschitz equivalence lies strictly between the two.
Indeed, it is possible for two finite polyhedra to be homeomorphic
but not bilipschitz equivalent, as observed in \cite{SS}, using the
examples from \cite{E, C1, C2}.  There are also the examples in
dimension $4$ now, as discussed in Section \ref{Homeomorphisms}.  In
the other direction there are finite polyhedra which are bilipschitz
equivalent but combinatorially distinct.  One can get this from
\cite{Si}, and this is clarified further by the general theory in
\cite{Su}.  Examples of the latter phenomenon could also be obtained from
\cite{Mi}, although for this it seems better to rely on the later
technology of $H$-cobordisms and $s$-cobordisms (as pointed out to me
by S. Cappell).  The analogous question for $4$-dimensional compact
smooth manifolds remains open, i.e., whether bilipschitz equivalence
implies smooth or piecewise-linear equivalence (which are the same in
this dimension).)

\section{Finite sets and geometry}

	Let us try think about geometry in a primitive way.  We are
given something like a metric space, or even just a set inside of some
${\bf R}^n$.  Maybe a curve or a surface.  Or maybe back to a finite
world, finite sets, or graphs.

	``Geometry'' is a tricky word.  One can make many different
kinds, different interpretations.  The idea of a metric space provides
one interpretation.  Distance and nothing more.

	Suppose that we are working inside some fixed ${\bf R}^n$.  We
can approximate any compact set with a finite set.  To make this
precise we can use the {\em Hausdorff metric}.  Given $A, B \subseteq
{\bf R}^n$, set
\begin{equation}
	D(A,B) = 
	   \max(\, \sup_{a \in A} \dist(a,B), \, \sup_{b \in B} \dist(b,A)).
\end{equation}
Thus $D(A,B) \le \epsilon$ means that every element of $A$ is within
distance $\epsilon$ from $B$ and vice-versa.  This defines a complete
metric on the set of nonempty compact subsets of ${\bf R}^n$.  If we
restrict ourselves to sets $A, B$ which are contained inside some
fixed compact set, then the resulting space is compact.

	Similar considerations apply inside any metric space.

	Inside this metric space of compact subsets the collection of
finite subsets is dense.  This is easy to show.  It provides a way to
make precise the idea that compact sets are simply ``limiting
versions'' of finite sets.

	We have considered compact sets as subsets of some larger
space, like a Euclidean space, but only for convenience.  Gromov has
given a version of the Hausdorff metric for compact metric spaces as
abstract spaces, independently of an embedding.  See
\cite{Gr0}.  Again metric spaces made from finite sets are dense
within the totality of all compact metric spaces.

	But is this really the way that we want to take limits
of finite sets?

	Let us work with subsets of ${\bf R}^n$.  Given a finite
subset $F$ of ${\bf R}^n$, let $\mu_F$ denote the probability measure
which is the sum of the point masses associated to elements of $F$
divided by their total number.  Thus $\mu_F(E)$ is simply the number
of elements of $F$ which lie in $E$, divided by the total number of
elements in $F$.  We may permit the elements of $F$ to have
repetitions.

	We can take limits of these measures to get other measures
on ${\bf R}^n$.  We shall use the following standard notion of
convergence.  If $\{\mu_j\}_j$ is a sequence of Borel measures on
${\bf R}^n$, each finite on bounded subsets of ${\bf R}^n$, and if
$\mu$ is another such measure, then we say that $\{\mu_j\}_j$
converges to $\mu$ if
\begin{equation}
	\lim_{j \to \infty} \int_{{\bf R}^n} \phi \, d\mu_j
		= \int_{{\bf R}^n} \phi \, d\mu
\end{equation}
for all continuous functions $\phi$ on ${\bf R}^n$ with compact
support.  A well known result in functional analysis implies that any
sequence $\{\mu_j\}_j$ of Borel probability measures -- nonnegative
measures with total mass $1$ -- has a subsequence which converges in
this sense.  The limiting measure may not be a probability measure,
though.  It will be nonnegative and have total mass at most $1$, but
the total mass could be strictly less than $1$.  Some mass might leak
out to infinity.  There are simple criteria to prevent this, such as
the requirement that the measures have support all contained in a
common compact set.

	This notion of convergence of probability measures leads us to
a different way to pass to the infinite realm with the notion of
finite sets.  It is not hard to show that any probability measure can
arise as the limit of measures associated to finite sets as above, and
in fact the same is true of subprobability measures (nonnegative Borel
measures with total mass at most $1$) if we allow portions of the
finite sets to leak off to infinity.

	Suppose now that we are given a sequence of finite sets $F_j$
in ${\bf R}^n$, and assume for simplicity that they are contained in a
fixed compact set.  We can form the corresponding measures
$\mu_{F_j}$, which we shall denote as $\mu_j$ for convenience.  Assume
that the $F_j$'s converge to a compact set $K$ in ${\bf R}^n$ in the
Hausdorff metric, and that the measures $\mu_j$ converge to a
probability measure $\mu$.  We can always achieve the existence of
these limits by passing to subsequences.  It is not hard to show that
the support of $\mu$ is contained in $K$.

	However, the support of $\mu$ may be a proper subset of $K$.
The measure $\mu$ encodes the locations of the points weighted
according to their multiplicity, and if the sets become too sparse in
some region, that part might disappear in the limit of the measures,
even if it persists in the limit of the sets.

	For instance, suppose that a point $p$ lies in each $F_j$.
Then $p$ will certainly lie in $K$ as well.  Suppose also that $p$ is
the only element of the ball $B(p,1)$ which lies in each $F_j$.  If we
allow points to have multiplicities, then let us ask that $p$ always
have multiplicity $1$ as well.  If the number of elements of $F_j$
tends to infinity as $j \to \infty$ -- the situation of primary
interest for us -- then $p$ will not lie in the support of the
limiting measure $\mu$, and in fact we shall have $\mu(B(p,1))=0$.
(For this it matters that $B(p,1)$ is the {\em open} ball with center
$p$ and radius $1$.)

	In some ways the idea of the limiting measure contains {\em
more} information than the idea of the limiting set, in some ways {\em
less}.  In any case they can be quite different, unless one imposes
additional conditions.

	Suppose now that we have sets in ${\bf R}^n$ which are not
finite but which are finite graphs, finite unions of line segments,
say.  We can take limits of them in the sense of Hausdorff limits, as
above, but maybe when we do that we lose some interesting geometric
information.

	Imagine for instance that we have a sequence of polygonal
curves $\Gamma_j$ which all connect some point $p$ to some point $q$.
We can choose the $\Gamma_j$'s so that they approximate the line
segment from $p$ to $q$ in the Hausdorff metric -- e.g., every element
of $\Gamma_j$ lies within a distance of $1/j$ of that line segment --
but so that the $\Gamma_j$'s are also pretty jagged, constantly going
up and down.  In such a way that each $\Gamma_j$ has twice the length
of the line segment which connects $p$ and $q$, for instance.

	Thus something is clearly lost in the limit, some finer
information about the way that the $\Gamma_j$'s behave.  In the limit
it is easy to go straight from $p$ to $q$, but not in any of the
approximations.  In fact, we could choose the $\Gamma_j$'s so that
they converge to the line segment in the Hausdorff distance but have
lengths tending to infinity.

	Here is another example.  Imagine a sequence of graphs in
${\bf R}^2$ which approximate a square in the Hausdorff distance.
That is already a bit funny, a sequence of one-dimensional sets converging
to a two-dimensional one, but it is easy to arrange.  In fact one can do it
in very different ways.  In one case one might approximate the square
by a bunch of parallel line segments.  They do not have to be connected
to each other.  In the limit ``connections'' emerge which are not at
all present in the approximations.  In the same way that a sequence of
finite sets can approximate a curve, or a square.

	On the other hand one might choose a sequence of
approximations which looks like a fine mesh.  A union of two
collections of parallel line segments, that is, vertical and
horizontal.  In this way all of the ``connections'' within the square
are approximated by ones in the graphs.  Approximated in the Hausdorff
metric; in terms of length the approximations may not work.  Suppose
that we are given two points in the square which are connected by a
line segment that is neither horizontal or vertical.  In order to go
between small neighborhoods of the two points in the approximating
graphs one would have to traverse a polygonal path that switches
between horizontal and vertical at least once.  The lengths of the
approximating paths would differ from the optimum by a definite
proportion.

	There is much more geometry in the approximations than we see
in the limit if we simply think in terms of Hausdorff limits or
limiting measures.  There are more elaborate theories of ``currents''
and ``varifolds'' which take into account additional information of
this type.  See \cite{Mo, Sm}.  In these theories a limit of something
like graphs will still be viewed as $1$-dimensional even if the limit
spreads out across a square or other higher-dimensional set.  These
theories take into account a notion of ``tangent directions'' in
addition to the distribution of the underlying points.  Both the
distributions of tangent directions and points are tracked through
measures, and they can become fuzzy in the limit, in the same way that
a sequence of finite sets can fill out arbitrary sets in the limit.

	Still these theories capture only some aspects of the
geometry.  Of course we cannot always expect to be able to extend a
given geometric idea to a limiting regime.

\section{Finitely generated groups}

	Suppose that we have a group $G$ with a finite number of
generators $g_1, \ldots, g_n$.  We can form {\em words} out of the
$g_i$'s, finite strings of $g_i$'s and their inverses, and all
elements of the group arise in this manner.  There may be many words
however which represent the same group element, and this can lead to
difficult combinatorial problems.  This is true even if the group is
finitely presented, so that the set of relations is generated by a
finite set of words.  The word problem may be algorithmically
unsolvable, for instance.  (See \cite{Ma}.)

	This is a situation with interesting combinatorics encoded in
finite constructions, the construction of words.  What are natural
limiting versions of these finite constructions?  One could simply
take infinite words, but that may not be so interesting after all. One
might hope for a natural kind of ``product integral'', at least in
some cases, but this is not clear.  For hyperbolic groups one has the
boundary at infinity as defined in \cite{Gr2}, and more generally one
can look at asymptotics of the group (in the word metric).  See
\cite{Gr3}.  There are also natural non-hyperbolic cases with good
``boundaries''.  See \cite{FM}.

	In the preceding examples we often had a continuous version
(or versions!) that we ``wanted'' to believe in, that are
well-established in mathematics.  In this case the matter is less
clear.

\section{Comments about symmetry and infinity}
\label{Comments}

	Let us accept the idea that there are vast regions of
unexplored territory in the world of finite and combinatorial
structures, and that transcendental/continuous mathematics is of
limited scope in comparison.

	It is often easier to prove theorems in the infinite world.
It is often easier to use the language more efficiently there to make
descriptions.  As in calculus, for instance; in addition to matching a
continuous idealization, it is simply easier to compute in calculus
than in discrete approximations.  Easier to make mathematical models
of physical processes, and to find clean formulations for both
questions and answers.

	One can also typically find more symmetries in mathematical
objects in the transcendental/continuous world.  Symmetries are
interesting in their own right and can facilitate the proving of
theorems.  We have discussed this before, in Chapter \ref{About
symmetry in language}.

	Without these kinds of symmetries it may be physically
impossible to find proofs, they would be just too long to find.  We
may also have let ourselves be trapped by the convenience of these
symmetries and continuous models.

\chapter{Measuring infinite complexity, I}

	How can we make measurements of different but infinite levels
of complexity?

	Ideas of infinite complexity have a prominent role in
classical analysis, even if they are not described explicitly as such.
Again we shall look at some specific examples.

\section{Functions and polynomials}

	Let us think about functions, on the real line for simplicity.
Real-valued functions.

	In learning calculus the first class of functions that one
might encounter are {\em polynomials}.  These are described by
finitely many parameters and can be evaluated at any given point
through finite arithmetic operations.  We might say that polynomials
have {\em finite complexity}.  Continuous or even smooth functions
in general have {\em infinite complexity}.  How can we be more
precise about this?  How can we make measurements?

	In calculus we learn about {\em Taylor polynomials} and {\em
power series expansions}.  These ideas provide ways to try to
approximate general functions by polynomials.  The sine, cosine, and
exponential functions have power series expansions which converge
everywhere and rapidly.  In some cases, like the tangent function, or
algebraic functions, there are power series expansions which work
locally but not globally, due to the presence of singularities.

	The same phenomenon occurs for rational functions, such
as $\frac{1}{1-x}$.  This blows up at $x = 1$ but for $|x| < 1$ it
is represented by the power series
\begin{equation}
	1 + x + x^2 + x^3 + \ldots
\end{equation}
One might argue that there is little point in using power series in
this case, though, since the function can already be ``computed''
through arithmetic operations directly.

	Another example in a similar vein is the function
$\frac{1}{i-x}$, where now we allow functions with complex values and
$i$ denotes the square root of $-1$.  This does not have a singularity
on the real line, but it does have one in the complex plane, at $x =
i$.  This means that the power series expansion at $x = 0$ cannot
converge for $|x| > 1$, because otherwise the function would not have
a pole at $x = i$.  One can check this directly anyway, but the point
is a general one.  We might try to approximate a function by
polynomials, through a power series expansion, and this can lead to
trouble even when there is not a singularity on the range being
considered.

	On the other hand we have the famous theorem of Weierstrass,
\index{Weierstrass' theorem}
to the effect that if $f(x)$ is any continuous function on a compact
interval $[a,b]$, then $f$ can be uniformly approximated by
polynomials.  That is, for each $\epsilon > 0$ there is a polynomial
$P(x)$ such that
\begin{equation}
	\sup_{x \in [a,b]} |f(x) - P(x)| < \epsilon.
\end{equation}
Of course the polynomial $P$ depends on $f$ and on $\epsilon$.
See \cite{Ru1} for a proof.

	If this works for any continuous function, why do we care so
much about power series?  Doesn't this already say that {\em all}
continuous functions are practically the same as polynomials?

	For the second question, the answer is yes in a sense and no
in a sense.  Continuous functions involve infinite processes in
general, and Weierstrass' theorem provides a way to say that they can
be approximated by finite processes.  Neither the {\em nature} nor
the {\em rate} of approximation provided by Weierstrass' theorem
are very strong however.

	To measure the rate of approximation we can proceed as
follows.  Given a positive integer $d$, let ${\cal P}_d$ denote the
collection of polynomials of degree at most $d$.  This is a vector
space of dimension $d+1$.  We can ask how close our given function
comes to this special subspace of the functions on $[a,b]$.  We make
measurements using the supremum norm, defined by
\begin{equation}
	\|g\| = \sup_{x \in [a,b]} |g(x)|
\end{equation}
when $g$ is a function on $[a,b]$.  For a given function $f(x)$
on $[a,b]$, consider the quantity
\begin{equation}
	\dist(f, {\cal P}_d) = \inf_{P \in {\cal P}_d} \|f - P\|
\end{equation}
defined for each $d$.  This measures the distance from $f$ to the
polynomials of degree at most $d$ in the supremum norm.  Weierstrass'
theorem says exactly that
\begin{equation}
	\lim_{d \to \infty} \dist(f, {\cal P}_d) = 0
\end{equation}
whenever $f$ is a continuous function on the compact interval $[a,b]$.
The converse is also true, because the uniform limit of continuous
functions is necessarily continuous.  (See \cite{Ru1}.)

	How fast do the numbers $\dist(f, {\cal P}_d)$ go to zero?

	The answer depends on the choice of function $f$.  The numbers
$\dist(f, {\cal P}_d)$ may go to zero as slowly as they like,
depending on the ``degree'' of the continuity of $f$.  This can be as
``bad'' as it likes.  To make this precise, recall that a continuous
function $f$ on a compact interval $[a,b]$ is always {\em uniformly
continuous}.  This means that for each $\epsilon > 0$ there is a
$\delta > 0$ so that if $x, y \in [a,b]$ and $|x - y| < \delta$, then
\begin{equation}
	|f(x) - f(y)| < \epsilon.
\end{equation}
For a general continuous function there is nothing that we can say about
the relationship between $\epsilon$ and $\delta$.

	Suppose that $f$ is everywhere differentiable on $(a,b)$.
Then for each $x, y \in [a,b]$ with $x \ne y$ there is a point $\xi$
in the open interval between $x$ and $y$ such that
\begin{equation}
	f(x) - f(y) = f'(\xi) \cdot (x-y).
\end{equation}
This is the {\em mean value theorem}.  From this we conclude that
\begin{equation}
	|f(x) - f(y)| \le \{\sup_{z \in (a,b)} |f'(z)| \} \cdot |x-y|
\end{equation}
for all $x, y \in [a,b]$.  The supremum on the right hand side
may not be finite, but it will be if $f'$ is continuous on $(a,b)$
and extends to a continuous function on $[a,b]$.  

	In this case we can take $\delta$ to be proportional to
$\epsilon$, i.e., we can take
\begin{equation}
	\delta = \epsilon \cdot (\sup_{z \in (a,b)} |f'(z)|)^{-1}.
\end{equation}
In general we might not be able to choose $\delta$ to be a linear
function of $\epsilon$.  For instance, if $[a,b] = [0,1]$ and $f(x) =
\sqrt{x}$, then we cannot take $\delta$ to simply be a multiple of
$\epsilon$ when $\epsilon$ is very small.  We have to take $\delta$ to
be at least $\sqrt{\epsilon}$, as one can check by taking $x = 0$.
For general continuous functions there is no universal law for the
relationship between $\epsilon$ and $\delta$.  Given any proposed
relationship one can find continuous functions for which it does not
work.

	One can bound the rate at which $\dist(f, {\cal P}_d)$ goes to
$0$ in terms of the ``modulus of continuity of $f$'', which means (in
essence) the way that $\delta$ depends on $\epsilon$.  If $f$ is
differentiable everywhere on $(a,b)$ and its derivative is uniformly
bounded, then the asymptotic rate will be faster than for $\sqrt{x}$,
and in the manner that one might expect (linear decay in the first
case, and not in the second).

	The relationship between the rate of decay of $\dist(f, {\cal
P}_d)$ and the modulus of continuity of $f$ is good but not perfect.
The modulus of continuity can never be better than linear, unless the
function is constant; for if $\delta$ were to go to $0$ faster than
any small multiple of $\epsilon$, then the derivative of the function
would have to vanish everywhere.  The rate of approximation of a
function by polynomials can be much faster than the reciprocal of the
degree however.  In effect it can take into account {\em higher}
derivatives as well.  For polynomials we simply have that $\dist(f,
{\cal P}_d) = 0$ for sufficiently large $d$.  If $f$ can be
represented by a power series that converges everywhere on the real
line, then we have very fast exponential decay, i.e.,
\begin{equation}
	\dist(f, {\cal P}_d) = O(\eta^d) \qquad\hbox{ as } d \to \infty
\end{equation}
for any $\eta > 0$, no matter how small.  If we take a function like
$f(x) = \frac{1}{x+i}$ (with its singularity off the real line), then
on any given interval $[a,b]$ we shall have exponential decay, but the
rate will depend on the interval.  One can verify this directly for
this particular function.  (As a general recipe one can employ power
series representations centered on disks in the complex plane, with
the center taken to lie off of the real line, at strategically-chosen
locations.)

	In general the rate of decay can be related to the size of the
derivatives of $f$, when they exist.  The more derivatives which exist
the better, the better for decay of $\dist(f, {\cal P}_d)$ that is.
One can also look for stronger forms of convergence, e.g., uniform
convergence of the functions and of their derivatives, or one could
look for approximations in other kinds of norms, integral norms for
instance.  These issues are thoroughly studied in analysis.  This is
part of {\em Approximation Theory}, about which many books have been
written.

	Another issue is the problem of being able to find
approximations effectively.  For the Weierstrass theorem and its
cousins there are fairly simple recipes, but there are other
situations where explicit solutions are not so easy to come by.  There
are nonconstructive methods in mathematics for proving the existence
of certain kinds of approximations, and this is an interesting
phenomenon for the broader themes of this book.  For Weierstrass'
theorem, for instance, we can rephrase the matter as saying that the
polynomials are {\em dense} in the (Banach) space $C([a,b])$ of
continuous functions on the interval $[a,b]$ equipped with the
supremum norm $\|\cdot\|$ defined above.  This provides a kind of
simplification in the language which permits other techniques to be
applied, even if it is less convenient for producing actual
approximations or measuring rates of decay.

	One line of reasoning would proceed as follows.  Suppose that
the polynomials were not dense in $C([a,b])$.  Then there would be a
{\em bounded linear functional} on $C([a,b])$ which vanishes on the
polynomials.  That is, there would be a linear mapping $\Lambda$ from
$C([a,b])$ into the scalars (real or complex numbers, as one prefers)
which is {\em bounded}, which means that there is a constant $M > 0$
such that
\begin{equation}
	|\Lambda(f)| \le M \cdot \|f\|
\end{equation}
for all $f \in C([a,b])$, and which satisfies
\begin{equation}
	\Lambda(P) = 0 	\qquad\hbox{ for all polynomials } P.
\end{equation}
The existence of such a $\Lambda$ is a standard consequence of the
{\em Hahn-Banach theorem}.  \index{Hahn-Banach theorem} (See
\cite{Ru2}.)  One then tries to show that this linear functional
cannot exist.  

	Note that the Hahn-Banach theorem relies on transfinite
induction, a very large and abstract kind of infinite process.

	The next step in this kind of reasoning is to say that
$\Lambda$ is actually represented by integration, i.e., there is a
finite measure $\mu$ (but not necessarily positive) on $[a,b]$ such
that
\begin{equation}
	\Lambda(f) = \int_a^b f \, d\mu \qquad\hbox{ for all functions }
							f \in C([a,b]).
\end{equation}
This uses a form of the {\em Riesz representation theorem}.  (See
\cite{Ru2}.)  

	In the context of Weierstrass' theorem this machinery is a bit
ridiculous, since there are much more direct ways to solve the
problem.  The ideas are illustrated better by the {\em M\"untz-Szasz}
theorem.  \index{M\"untz-Szasz theorem} Suppose that we are working on
the unit interval $[0,1]$, and we are given a set $E$ of positive
integers.  Let ${\cal P}(E)$ denote the set of polynomials which are
linear combinations of constant functions and monomials $x^j$, $j \in
E$.  Under what conditions on $E$ is ${\cal P}(E)$ dense in $C([0,1])$
dense?  The answer is that this happens if and only if
\begin{equation}
	\sum_{j \in E} \frac{1}{j} = \infty.
\end{equation}
This is the {\em M\"untz-Szasz theorem}.  See \cite{Ru2} for a proof,
which follows the lines above and then uses complex analysis to analyze
the relevant measures.

	Nonconstructive methods are less satisfactory for making
measurements of the level of infinite complexity, but there are
circumstances in which they can be used to analyze the rate of
approximation in a useful way.

\section{Polynomials, schmolynomials}

	Do we really want to measure the complexity of functions in
terms of approximation by polynomials?  Polynomials provide a natural
notion of ``finite complexity'', but not the only one.  For this
discussion we restrict ourselves to functions defined on the real
line, say, or intervals within the real line.

	Instead of using polynomials we could take rational functions
(with poles off the domain) as our basic ``finite'' objects, as we
mentioned earlier.  We can try to measure the complexity of functions
in terms of their approximations by rational functions.

	We want to begin with some building blocks.  We start with the
rational function $\frac{1}{1+x^2}$, and then we add two parameters as
follows.  Given real numbers $a$ and $b$, consider the function (in $x$)
defined by
\begin{equation}
\label{bump}
	\frac{b^2}{b^2 + (x-a)^2}.
\end{equation}
This function looks like a kind of bump.  It is easy to see that this
function assumes its maximum at $x = a$ and steadily decreases as $x$
moves away from $a$ in either direction.  The value of this function
at $x = a$ is $1$.  As one moves away from $a$ the function starts to
become small.  The rate at which it does this is measured by $b$.
Indeed, notice that our bump function is roughly constant at the scale
of $b$; given any interval of length $b$, the function does not change
very much.  On intervals of length larger than $b$ it does change more
significantly, the values steadily going down to zero with each step
away from $a$.

	The rate of decay is not real fast, but still pretty fast.
One can make it faster by using higher powers than $2$.  One might
hope for an exponential rate of decay, but of course one cannot
achieve that with rational functions.  One could use exponentials
instead, like Gaussians.

	At any rate our bump function (\ref{bump}) provides a very simple
method for approximating continuous functions by rational functions.
Fix a $b > 0$, very small, and think of approximating a given function
$f(x)$ on the unit interval $[0,1]$ (say) in the following manner.
Assume for simplicity that $b$ is the reciprocal of a large integer,
$b = N^{-1}$.  We can approximate $f$ by a function of the form
\begin{equation}
	\sum_{j=1}^N 	\phi_j \cdot 
		\frac{\frac{j}{N}^2}{\frac{j}{N}^2 + (x-\frac{j}{N})^2}
\end{equation}
for some choice of numbers $\phi_j$.  This is reasonably intelligent
precisely because our bump function (\ref{bump}) looks like a bump, it
is roughly concentrated in the interval $[a-b, a+b]$, and it takes
roughly the value $1$ there.  In our case we are taking the $a$'s to
be multiples of $b = N^{-1}$, and this means that our bump functions
\begin{equation}
\label{sum of bumps}
	\frac{\frac{j}{N}^2}{\frac{j}{N}^2 + (x-\frac{j}{N})^2}
\end{equation}
should be roughly independent of each other.  That is, they should
not interact too much, because they decay as one moves away from
their centers.

	How should we choose the coefficients $\phi_j$?  An obvious
guess might be to take $\phi_j = f(\frac{j}{N})$.  This does not work
so well, even for constant functions.  There are some multiplicities
that we have to take into account, but that is the only real problem.
To see this define a ``multiplicity function'' $M(x)$ by
\begin{equation}
	M(x) = \sum_{j=1}^N 	
		\frac{\frac{j}{N}^2}{\frac{j}{N}^2 + (x-\frac{j}{N})^2}.
\end{equation}
This is what we would get in the sum by taking all the coefficients
$\phi_j$ to be $1$.  For our given function $f$ let us choose the
coefficients to be given by
\begin{equation}
	\phi_j = f(\frac{j}{N}) \cdot M(\frac{j}{N})^{-1}.
\end{equation}
Denote by $f_N$ the function that results from (\ref{sum of bumps})
with these coefficients.  Then we have that $f_N$ converges to
$f$ uniformly on $[0,1]$ as $N \to \infty$.  Let us see why this is
so.

	Fix a point $x$, and consider $f_N(x)$, which is given
explicitly by
\begin{equation}
	f_N(x) = \sum_{j=1}^N 	f(\frac{j}{N}) \cdot M(\frac{j}{N})^{-1} \cdot
		\frac{\frac{j}{N}^2}{\frac{j}{N}^2 + (x-\frac{j}{N})^2}.
\end{equation}
This is a weighted average of the values $f(\frac{j}{N})$ of $f$.  It
is an {\em average} because the weights are positive and their sum is
$1$, the latter a consequence of the presence of the multiplicity
function.  The point now is that if $N$ is very large then these
weights are concentrated at the $j$'s for which $x-\frac{j}{N}$ is
small compared to $N^{-1}$.  When $x-\frac{j}{N}$ becomes large
compared to $N^{-1}$, the weights become small, and the contribution
of those terms is small.  In short $f_N(x)$ is an average of values of
$f$ with the average concentrated at values of $f$ taken at points
$\frac{j}{N}$ which are very close to $x$ when $N$ is large enough.
One can then use the (uniform) continuity of $f$ to say that this
average is (uniformly) close to $f(x)$.

	One needs to make some computations to make this argument
precise, but this is not very difficult.  One can also see clear
relationships between the modulus of continuity of the given function
and the rate of approximation of the $f_N$'s.  These can be formulated
as statements about the relationship between the modulus of continuity
and the rate of approximation as measured by the {\em degree} of the
rational functions (in this case the $f_N$'s).  The main point is that
the parameter $b$ (in this case $1/N$) corresponds to the resolution
of the approximation in a way that can easily be made precise.

	This is one way to make approximations by rational functions,
but it is by no means the only way, or the most efficient way.  There
are better methods for taking into account more smoothness or more
subtle oscillations, but we shall not pursue this here.  One of the
main points is to use rational functions and parameters like $a$ and
$b$ to treat {\em locations} and {\em scales} simultaneously and
individually.

	Still this gives a flavor, and a way to look at measuring
complexity of functions.  Note that some functions admit much better
approximations by rational functions than by polynomials.  Rational
functions can deal with singularities more efficiently, and our bump
functions (\ref{bump}) are localized much better than polynomials ever
would be.

	If one likes localization, one might as well ask about
approximation by step functions instead.  Recall that a {\em step
function} is simply a linear combination of characteristic functions
of intervals.  One might also describe them as being
piecewise-constant.  It is very easy to approximate a continuous
function on a compact interval uniformly by step functions, simply
using the uniform continuity of the function.  Step functions are also
described by only a finite number of parameters, and so provide
another notion of ``finite complexity'' from which to work.

	Of course step functions are also discontinuous, which causes
trouble in some contexts.  One might prefer to approximate functions
by others which are at least as regular as the original.  Thus one
might consider instead approximations by {\em piecewise linear
functions}.  For these it is reasonable to demand continuity when
approximating continuous functions.  For piecewise-linear functions
one gives up differentiability at the break points.  To get higher
degrees of smoothness one can look at functions which are piecewise
polynomial, like {\em polynomial splines}.

	Another deficiency of approximation by step functions is that
the rate of approximation will never be so great, never more than
linear basically.  By using piecewise linear functions or splines of
higher order one can get faster rates of approximation of functions
which have some smoothness.

	The idea that we want to approximate functions by functions
which are as good as the original points to a deficiency in the
earlier method for approximating functions by rational functions.
While it enjoyed a certain naive simplicity, it suffers from the fact
that to get a good approximation we use bump functions which are
highly concentrated, with the resolution parameter $b$ chosen very
small, so that our basic building blocks are relatively singular.  They
have large derivatives, for instance.

	There are more efficient methods of approximation for avoiding
this problem.  The area of harmonic analysis often called
``Littlewood-Paley theory'' deals with this issue in a good way.
Basically one should make approximations that are adapted to scale and
location, and not just individual locations.  Thus if a function is almost
constant in some region, it is not very efficient to approximate it
by highly concentrated bumps there.  Now-a-days these matters are
typically discussed in the context of ``wavelets'' ({\em ondelettes}
in French).

	One can make other methods of approximation.  A particularly
important one is to use {\em Fourier series}.  Suppose that we have a
function on the real line, continuous say, which for simplicity is
periodic with period $1$.  One can try to realize it as a sum of
simpler functions of the form $\sin 2 \pi n x$, $\cos 2 \pi n x$,
where $n$ is a positive integer.  One should also include constant
functions.

	There are many results about being able to do this.  Let us
call a function a ``trigonometric polynomial'' if it is a finite
linear combination of sines, cosines, and constants as above.  A
variant of the Weierstrass approximation theorem states that
continuous periodic functions can always be approximated uniformly by
trigonometric polynomials.  If one works instead in $L^2$ one gets an
orthogonal basis from the trigonometric functions, so that an $L^2$
function can always be expressed as an infinite series of
trigonometric functions with the coefficients computed by simple
integrals.  The series will converge in the $L^2$ norm, but it may not
converge uniformly even if the original function were continuous.
These are well-studied issues, with many results known.  See
\cite{Ru2, Z}, for instance.
						
	In a certain sense the analysis of functions by trigonometric
functions is very similar to analysis through polynomial approximation.
One can view periodic functions as being functions on the unit circle
in the plane ${\bf R}^2$, and trigonometric polynomials on the
line correspond exactly to the usual polynomials on ${\bf R}^2$
restricted to the unit circle.

	In other ways these various kinds of approximation and
expansions are very different.  One of the reasons that Fourier series
are very useful is that partial differential equations, such as the
equation which models the vibration of a string, become much simpler
when one expands functions in terms of Fourier series.  An individual
sine or cosine function can represent a ``standing wave'', and the
frequency $n$ is very important, corresponding to pitch for sound.
These things do not work as nicely for ordinary polynomials, or for
rational functions.

	On the other hand Fourier series are not very good for
treating local behavior in an efficient manner.  If one has a function
that is very calm most of the time but then erratic at a particular
location, then it is hard to see that in the Fourier series, because
the effect of the erratic behavior spreads out.  By using rational
functions, or splines, or wavelets, one can treat local phenomena
differently.

	One might also have a large region in which the function
is very calm, another where it is not.  Wavelets are pretty good at
figuring out which is which, and treating them automatically.  That is,
they are good for accommodating calm or erratic behavior that depends
on the function in an automatic way.

	There is a more elaborate way to try to ``analyze'' the
behavior of functions, based not on a (linearly independent) basis of
basic functions, but on a larger collection from which one can make
a library of bases.  Different bases to see different kinds of behavior,
be it waveform, localized oscillations, something more spread out.
One then has algorithms for choosing the best basis from the library
(where ``best'' has to be explained).  See \cite{Co, CMQW, CW1, CW2}.

	In short, there are many different methods for ``analyzing''
functions, which lead to different ways to measure complexity.

\section{Linear operators}

	Let us consider a very concrete linear operator which arises
in mathematics.  Let $f(x) = f(x_1,x_2)$ be a function on ${\bf R}^2$,
and consider the operator defined by
\begin{equation}
\label{definition of R}
	R f(x) = {\rm p. v.} \int_{{\bf R}^2} f(y) \, 
			\frac{(x_1-y_1)(x_2 - y_2)}{|x-y|^4} \, dy.
\end{equation}
This is a {\em singular integral operator}, because the integral does
not exist in the usual sense.  One has to define it as a {\em
principle value} (as in Section \ref{The Hilbert transform}).  This
means that in fact we set
\begin{equation}
	R f(x) = \lim_{\epsilon \to 0} \int_{|x-y| > \epsilon} f(y) \, 
			\frac{(x_1-y_1)(x_2 - y_2)}{|x-y|^4} \, dy.
\end{equation}
If $f$ is smooth and compactly supported (or even much less) one can show
that this limit always exists, by rewriting it as
\begin{eqnarray}
\label{adjusted}
	&& \int_{|x-y| > 1} f(y) \, 
			\frac{(x_1-y_1)(x_2 - y_2)}{|x-y|^4} \, dy  +  	\\
  && \lim_{\epsilon \to 0} \int_{\epsilon < |x-y| \le 1} (f(y) - f(x)) \,
			\frac{(x_1-y_1)(x_2 - y_2)}{|x-y|^4} \, dy.
								\nonumber
\end{eqnarray}
To make this rewriting we use the fact that
\begin{equation}
	\int_{\epsilon < |x-y| \le 1}
		\frac{(x_1-y_1)(x_2 - y_2)}{|x-y|^4} \, dy = 0,
\end{equation}
which can be checked directly.  (It is easiest to see this using symmetry
considerations, noticing how the integrand changes sign when one reflects
across an axis.)  The point now is that if $f$ is smooth and compactly
supported, then 
\begin{equation}
\label{estimate}
	\Bigl| (f(y) - f(x)) \, \frac{(x_1-y_1)(x_2 - y_2)}{|x-y|^4} \Bigr|
		\le \{\sup_{z \in {\bf R}^2} |\nabla f(z)| \} \cdot |x-y|^{-1}.
\end{equation}
This follows from the fact that $|f(y) - f(x)|$ is no greater than
$|x-y|$ times the supremum of $|\nabla f|$ on the line segment which
joins $x$ to $y$.  Once we have this estimate (\ref{estimate})
it is easy to show that the limit in (\ref{adjusted}) exists.
(Indeed now the ``singular'' integral converges even if we put absolute
values inside the integral, which was not true before.)

	Thus the singular integral exists for nice functions $f$.
This operator arises mathematically in the following way.  Recall that
the Laplacian $\Delta$ is the second-order differential operator
defined by
\begin{equation}
	\Delta = \frac{\partial^2}{\partial x_1^2} 
			+ \frac{\partial^2}{\partial x_2^2}
\end{equation}
(since we are in two dimensions).  Suppose that we have a smooth
function $f$ with compact support, and we want to be able to derive
information about the various second derivatives of $f$ from the
knowledge of $\Delta f$ alone.  A basic formula is
\begin{equation}
\label{formula}
	\frac{\partial^2}{\partial x_1 \partial x_2} f
		= R(\Delta f)
\end{equation}
(except for a multiplicative constant).  Thus $R$ is really the
quotient of $\frac{\partial^2}{\partial x_1 \partial x_2}$ and
$\Delta$.  This formula holds for practically any function, 
modulo suitable interpretation.

	This kind of identity can be obtained as in \cite{St1}, as
well as many others like it.  For the moment the main point is that
this is a linear operator of basic importance in mathematics, and
there are many others like it.  What can we say about it in terms of
complexity, and the infinite processes which it involves?  Is it
``simpler'' than an arbitrary operator, can it be described 
more ``efficiently'' than in (\ref{definition of R}) above?

	In this case the answer is yes.  Our operator enjoys extra
symmetries.  It commutes with translations, for instance.  This
implies that it can be ``diagonalized'' by the Fourier transform.
More precisely, if $f(x)$ is a function on ${\bf R}^2$, then its {\em
Fourier transform} is the function $\widehat{f}(\xi)$ on ${\bf R}^2$
defined by
\begin{equation}
	\widehat{f}(\xi) = 
	\int_{\bf R} f(x) \, e^{- 2 \pi i \, \langle x \, \xi \rangle} \, dx.
\end{equation}
We are using here the standard inner product $\langle \cdot \, \cdot
\rangle$ on ${\bf R}^2$.  This is like the definition given in Section
\ref{The Riemann-Lebesgue lemma} for functions on the line, and the inversion
of the Fourier transform works in the same manner as indicated there.
For our operator $R$ we have that
\begin{equation}
	{R(f)}^{\widehat{}}\, (\xi) = 
			\frac{\xi_1 \xi_2}{|\xi|^2} \, \widehat{f}(\xi),
\end{equation}
except for a muliplicative constant.  In other words, the complex
exponentials are eigenfunctions for $R$, and they are a ``complete
set''.  (We are ignoring some technical details here, but they are
standard.  See \cite{SW} for more information about Fourier transforms
and this fact that they diagonalize operators which commute with
translations.)

	The fact that we have a diagonalization is a way to say that
there is less ``complexity'' in this matrix than we might have
thought.  A general $N \times N$ matrix has $N^2$ entries, a diagonal
matrix has $N$ entries.  One could make suitable discretizations in
the present story to get closer to matrices if one wanted to.

	In fact all sorts of operators are ``diagonalizable'' , in the
same way that real symmetric matrices can be diagonalized, by the
spectral theory.  (See \cite{Ru3}, for instance.)  This sounds pretty
good, but for considerations of complexity one should also take into
account the transformation (``change of basis'') needed to effect the
diagonalization.  In this respect the Fourier transform behaves fairly
well.  It can be computed efficiently, through the fast Fourier
transform.

	Thus one might say that the complexity of our operator $R$ is
much less than it might appear to be at first, or than it would be for
an arbitrary operator.  There are still some problems with this,
however.  The Fourier transform works very well for some kinds of
information but not others.  It is not easy to see localizations in
the Fourier transform.  By contrast it turns out that $R$ does behave
fairly well in terms of localization.  That is, the effect of $R$ on a
function in a certain region is not too sensitive to what happens far
away.  The Fourier transform is not like that, it is more global.

	The precise localization properties of $R$ are a bit tricky.
Some of them are described in \cite{St1}.  They are reflected in part
by the boundedness of $R$ on many natural Banach spaces of functions.
This is not at all clear from the diagonalization through the Fourier
transform, because the action of the Fourier transform on most
function spaces is not so easy to understand.

	For example, $R$ defines a bounded linear operator on
$L^p({\bf R}^2)$ when $1 < p < \infty$.  This means that for each such
$p$ there is a constant $C(p)$ such that
\begin{equation}
	\int_{{\bf R}^2} |R f(x)|^p \, dx 
		\le C(p) \int_{{\bf R}^2} |f(x)|^p \, dx 
\end{equation}
for all functions $f$ on ${\bf R}^2$.  This kind of result is very
useful in analysis, especially partial differential equations; it
provides a way to control the size of the $\frac{\partial^2}{\partial
x_1 \partial x_2}$ derivative of a function in terms of the size of
the Laplacian.  It does not work for $p = 1, \infty$, and in
particular it is better to measure ``size'' here in terms of certain
types of averages instead of the supremum.  There are many other
natural function spaces on which $R$ is also bounded.

	For $p = 2$ this boundedness can be proved using the Fourier
transform.  The main points are that the Fourier transform preserves
the $L^2$ norm, and on the Fourier transform side $R$ looks like the
operator of multiplication by a bounded function.  This does not work
for other $L^p$ spaces because the Fourier transform does not
have such simple behavior on them.  The fact that $R$ does behave
well on them is a reflection of its sensitivity to localization
and the underlying geometry of ${\bf R}^2$.

	There are other methods from Harmonic Analysis for analyzing
operators which avoid some of the difficulties suffered by the Fourier
transform.  Now-a-days people often speak in terms of ``wavelets'',
which provide special bases for function spaces which have moderate
sensitivity to both location and oscillation while at the same time
having natural scale-invariance.  For instance, one can characterize
$L^p$ spaces in terms of wavelets in a reasonably simple way when $1 <
p < \infty$, purely in terms of the {\em sizes} of coefficients.  This
is not true for the Fourier transform when $p \ne 2$.  This is called
the {\em unconditional basis property}, and it is a good substitute
for orthogonality which is meaningful in non-Hilbert spaces.

	It turns out also that when we represent $R$ as a matrix in
terms of a wavelet basis it is almost diagonal in a precise sense
which is useful in practice.  

	In fact wavelets cooperate well with a variety of linear
operators and function spaces.  The Fourier transform is better for
getting very sharp information in very special situations, while
wavelets are typically less precise but more flexible.  Ideas related
to wavelets are often applicable in cases where there is no
translation-invariance or even a group structure necessarily.  We have
emphasized the operator $R$ for the sake of illustration, but the
methods are actually very general.

	See \cite{Db2, Db3, M1, M2} for more information about
wavelets and the analysis of functions and operators.  Note that they
are powerful tools simply for studying the complexity of functions,
independently of operators.

	Many of these ideas for analyzing the complexity of functions
and operators existed in Harmonic Analysis long before the more recent
advent of wavelets.  One of the main points about wavelets is that one
can get actual {\em orthogonal} bases with these good properties;
before one had various methods for breaking functions and operators
into simple pieces and putting them back together again, and these
methods enjoyed much the same properties in terms of localization in
space and frequency, but the orthogonality properties were not as
exact.  Orthogonality is often very convenient even in situations
where its exact form is not crucial.

	Wavelets also work well for numerical computation.  In part
this comes from their orthogonality, which eliminates some of the
redundancies needed as well as providing exact formulas that do not
require additional analysis.  Another important point is that the
wavelets themselves can be computed efficiently, and the
representations of functions in terms of wavelets can be computed
efficiently.  This uses the particular structure of the wavelets
defined by Daubechies in \cite{Db1}.  It is not that there is a simple
formula for these wavelets, but rather an iterative construction which
lends itself to numerical computation.  For applications to the
efficient computation of operators like $R$, see \cite{BCR, Co}.

\section{General comments about operators}

	Let us look a bit more about general issues of ``infinite
complexity'' related to linear operators which arise naturally in
mathematics.

	Notice that an operator like $R$ can be discretized in
``meaning'' as well as in the formula for its definition.  One could
work on a discrete $2$-dimensional lattice, define discrete versions
of the Laplacian and the mixed second derivatives, and look for a
``quotient'' as before (as in (\ref{formula})).  One would still get
something which looks approximately like $R$.

	Thus in a sense the concept of $R$ is not unavoidably
infinite.  The ``infinite'' realm is more convenient for certain types
of theoretical considerations, though.

	There are other aspects of operator theory that do require the
infinite world in a stronger way.  For this discussion I wish to stick
to fairly ``concrete'' operators, of the sort which arise in
connection with partial differential equations, as opposed to more
abstract operator theory.

	In the ``pseudodifferential calculus'' one has a collection of
linear operators which includes standard differential operators but
which also allows more complicated constructions.  One has a way to
try to take inverses, square roots, and other operations of interest.

	Actually, one does not perform these operations in an {\em
exact} way, but instead one has a symbolic calculus for converting
questions about operators into simpler questions about functions, at
the expense of allowing certain kinds of errors.  The composition of
two operators corresponds to the multiplication of a pair of functions
(the ``principal symbols'' of the operators), but one can recover the
operators from the principal symbols only to first approximation,
modulo an error which is of ``lower order''.  One can make more
refined calculations to make the errors even more modest, but the
improvements are in terms of {\em smoothness} rather than {\em size}.

	In other words, by allowing certain types of errors one can make
symbolic calculations which detect much of the mathematically significant
information in a convenient manner.  One has a good theory here of the
``infinite'' part of the complexity.

	For theoretical purposes this can be quite powerful.  One
might be interested in tracking singularities, for instance, or one
might not mind being able to invert an operator only up to a finite
rank error.  In some cases such ambiguities are unavoidable.

	For the matter of smooth or finite rank errors we see a very
clear difference between infinite mathematics and finite
approximations.  For infinite mathematics finite rank errors are not
so bad, but in a discrete approximation or numerical computation of
actual solutions the matter is far different.  See \cite{BCR2} for a
treatment of pseudodifferential calculus in a computational setting.

	Similarly the distinction between {\em smoothness} and {\em
singularity} reflects an important difference between finite and
infinite situations.  In the mathematical study of partial
differential equations which describe the evolution of some object,
like a fluid, one is very interested in knowing about the possibility
and nature of ``blow-up''.  In finite realms the distinction is much
less sharp.

	For the question of invertibility, one of the basic points
behind this discussion is the following.  Suppose that we have an
operator $T$ acting on some Banach space, and that we know that it is
a perturbation of a simpler operator $T_0$ that we understand fairly
well already.  Let us write $E$ for the error, so that $T = T_0 + E$.
If we know that $T_0$ is invertible, then under what conditions on the
error term $E$ can we invert $T$?

	This is true if $E$ is sufficiently small (in operator norm),
by a well-known computation.  In the context of differential
equations, or the pseudodifferential calculus, one typically knows
instead that the error term is somehow {\em smooth} instead of small.
In practice this leads to some form of {\em compactness}, like the
statement that $E$ can be approximated by finite rank operators in the
operator norm.  Thus we can proceed to $T = T_1 + E_1$, where $T_1$ is
a small perturbation of $T_0$ and therefore invertible, while $E_1$ is
now of finite rank.

	Given the knowledge that $T$ is a finite rank perturbation of
an invertible operator, in order to show that it is invertible it
suffices to show that it is either injective or surjective.  This is
true automatically in finite dimensions but not in general in infinite
dimensions.

	In the context of partial differential equations this fact can
be very useful.  The problem of inverting an operator might correspond
to the problem of solving a differential equation with some boundary
conditions, for example.  The preceding facts can permit one to reduce
a question of existence to one about uniqueness, for instance, which
might be approachable by direct calculation.

	Thus in the end the mathematical theory sometimes clears away
the brush until one comes down to a finite problem.  The mathematical
theory often stops then.

	There is another aspect of ``infinite complexity'' in operator
theory, in which one has an operator in hand with some compactness
properties and one wants to analyze more precisely the way that it can
be approximated by finite-rank pieces.  This issue is analogous to the
one of approximating functions by simpler pieces, and indeed there are
various ways in which the two questions are closely connected.  See
\cite{So, RS}, for instance.

\chapter{Measuring infinite complexity, II}

	We want to continue looking at mathematical tools for measuring
infinite complexity, but now with more emphasis on geometry and topology.

\section{Sizes of sets}
\label{Sizes of sets}

	In the finite realm one might simply measure size by counting
the number of points, but for infinite sets the matter is more
complicated.  We certainly want to say that a plane is bigger than a
line even if they are the same in terms of cardinality.

	Consider a subset $E$ of some ${\bf R}^n$.  How can we measure
its size?  The following is a simple way.  For each $\epsilon > 0$ let
$N(\epsilon)$ denote the smallest number of balls of radius $\epsilon$
that it takes to cover $E$.  This will be finite as soon as $E$ is
bounded.  One then looks at the asymptotic behavior of $N(\epsilon)$
as $\epsilon \to 0$.  For instance, if
\begin{equation}
\label{minkowski}
	N(\epsilon) = O(\epsilon^{-d}) \quad\hbox{ as } \epsilon \to 0
\end{equation}
for some $d \ge 0$, then this is a way to say that $E$ is at most
$d$-dimensional in terms of size.

	For example, a line segment satisfies (\ref{minkowski}) with
$d = 1$, and a square satisfies (\ref{minkowski}) with $d = 2$.  For
an entire line $N(\epsilon)$ is always infinite, as it is for any
unbounded set.  In that case it is better to restrict oneself to
bounded subsets and make measurements there, and then look to see how
the bounds depend on the subset, e.g., on its diameter.

	In this kind of measurement the ``dimension'' parameter $d$
need not be an integer.  For instance, for the usual Cantor set the
correct $d$ is $\log 2 / \log 3$.

	Of course one is not restricted to taking powers in
(\ref{minkowski}), one can use more complicated functions too.

	This method of measuring size has some disadvantages.  Suppose
that $E$ is a countable set whose closure contains a ball.  Then the
best that one can do with $N(\epsilon)$ is say that it satisfies
(\ref{minkowski}) with $d = n$, which is true of any bounded subset of
${\bf R}^n$.  This is not so nice, and one might want to say that such
sets are very small.  For this one has to look more deeply into the
structure of the set.  Let us describe a couple of ways for doing
this.

	Given $d \ge 0$, set
\begin{eqnarray}
	A^d(E) = \inf\{ \sum_j (\diam X_j)^d : \{X_j\}_j 
					\hbox{ is a sequence of subsets}
									\\
		\hbox{of } {\bf R}^n \hbox{ which covers } E\}.
								\nonumber
\end{eqnarray}
This is called the $d$-dimensional {\em Hausdorff content} of $E$.
\index{Hausdorff content}
It differs from the above because one can look at coverings of
variable size, i.e., the $X_j$'s do not have to all have the same
diameter.  For instance, $A^d(E) = 0$ whenever $E$ is countable and $d
> 0$.  For something like a line segment or a square the Hausdorff
content would coincide with the usual measure, modulo a normalizing
factor, and assuming that one took the right dimension $d$.

	A problem with the Hausdorff content is that it is always
finite, at least for bounded sets.  A square has finite $1$-dimensional
content, for instance; the finiteness of $A^d(E)$ does not really
say something about $E$, in terms of the dimension $d$.  By contrast,
a condition like (\ref{minkowski}) does see the dimension in a strong
way, (\ref{minkowski}) is not true for a square with $d=1$.

	The Hausdorff content can see the dimension in a different
way.  If $S$ is a square then $A^d(S)$ is always finite, but it is $0$
as soon as $d > 2$.  More generally, if a set $E$ satisfies
(\ref{minkowski}) for some $d$, then one has $A^s(E) = 0$ for all $s >
d$.

	In the Hausdorff content the $X_j$'s do not have to be small.
By contrast the asymptotics of $N(\epsilon)$ as $\epsilon \to 0$ must
unavoidably reflect the behavior of the set at small scales.  We can
do this also for Hausdorff content through a more precise definition.

	Let $d \ge 0$ and $r > 0$ be given.  Again $d$ corresponds
to dimension, while $r$ will correspond to a length scale.  Set
\begin{eqnarray}
H^d_r(E) = \inf\{ \sum_j (\diam X_j)^d : \{X_j\}_j \hbox{ covers } E
					\hbox{ and }
									\\
				\diam X_j < r \hbox{ for all } j\}.
								\nonumber
\end{eqnarray}
In other words, $H^d_r(E)$ is the same as $A^d(E)$ except that
all the $X_j$'s are required to be small.  It is easy to see that
$H^d_r(E)$ increases as $r$ decreases, because one is then
taking infima over smaller sets.  Set
\begin{equation}
	H^d(E) = \lim_{r \to 0} H^d_r(E).
\end{equation}
This defines the $d$-dimensional {\em Hausdorff measure} of $E$.
\index{Hausdorff measure}
The limit always exists because of the monotonicity in $r$,
but we allow the value $\infty$.

	This is now more sensitive to the dimension.  If $S$ is a
square then $H^2(S)$ is the same as the usual $2$-dimensional measure
of the square, except perhaps for a normalizing constant, but $H^d(S)
= \infty$ when $d < 2$, as it should be.  Similar remarks hold for
line segments or the standard Cantor set, with respect to the correct
dimension in each case.

	If $E$ satisfies (\ref{minkowski}) for some $d$, then $H^d(E)
< \infty$ for that choice of $d$.  The converse is not true.  Countable
sets have $0$ $H^d$-measure for each $d > 0$, but they can be more
complicated for (\ref{minkowski}).  One does not even have to take
countable dense sets, simply a sequence covering to a point will
work if the convergence is not too fast.  Similarly one can look
at sequences of line segments, or little squares.

	Indeed, one of the main advantages of Hausdorff measure and
Hausdorff content over measurements as in (\ref{minkowski}) is that
they are {\em countably subadditive}; the size of a countable union is
less than or equal to the sum of the measures.  One of the main
advantages of Hausdorff measure over Hausdorff content is that
Hausdorff measure is {\em countably additive}, the measure of a
countable disjoint union of {\em measurable} sets is the sum of the
individual measures.  The failure of countable (or even finite)
additivity of Hausdorff content reflects the fact that it is not
required to see small scales.  When $d=0$, for instance, $H^d$ simply
counts the number of elements in the given set, while $A^d$ can be
much smaller.

	On the other hand countable additivity does allow ``messy''
sets with many little pieces floating around.  In some situations one
might want to forbid that, and measurements in terms of $N(\epsilon)$
can be more appropriate.  It depends on the situation.

	Notice that $H^d(E) = 0$ if and only if $A^d(E) = 0$.  We have
\begin{equation}
	A^d(E) \le H^d(E)
\end{equation}
for all sets, just from the definitions, but the reverse inequality is
wrong, as we have seen.  However, $H^d(E) = 0$ when $A^d(E) = 0$,
because in that case the efficient coverings of $E$ for computing
$A^d(E)$ had to involve sets of small diameter.

	The {\em Hausdorff dimension} \index{Hausdorff dimension} of
$E$ is defined by
\begin{equation}
	\inf \{ d \ge 0 : H^d(E) = 0\}
\end{equation}
Note that $H^d(E) = 0$ for all $d > n$ when $E$ is contained in ${\bf
R}^n$, and so the Hausdorff dimension is no greater than $n$ in this
case, which is good.  One can also talk about the {\em (upper)
Minkowski dimension} as the smallest $d$ for which (\ref{minkowski})
holds, although for unbounded sets it is better to look at the
smallest $d$ such that (\ref{minkowski}) holds for all bounded
subsets.  There are other notions of dimensions and other measures
related to these, as in \cite{Fa, Fe, Mat}.

	Sets in Euclidean spaces can be complicated.  Very
complicated.  We have many ways to measure how complicated they are.
So far we have talked about size, but indeed also the geometry in a
slightly more subtle way, through the structure of the coverings, the
different ways that coverings can be organized.  These coverings speak
about the geometry of a set, maybe in a tricky way if we allow
coverings by sets of different diameters, as in Hausdorff measure and
Hausdorff content.  

	These ideas can make sense in ``finite'' contexts, in
discretizations, but it is not completely easy to do something beyond
simply computing the number of elements in a set.  One should
distinguish between scales in the discretization, e.g., by working
with a version of $H^d_r$ in which $r$ is much larger than the
smallest scales of the discretization but still very small.
Otherwise, if one simply takes $r$ to be the minimal scale of the
discretization, one cannot do more than simply count points.

	Even this can be useful, especially if one makes additional
localizations, e.g., counting the number of points in the set inside
balls of certain radius.  Looking at this as a function of the radius,
comparing the result with its counterpart for the whole space, etc.

	Still it is much easier to distinguish between different
scales in infinite mathematics than in finite worlds. 

	Even for questions of size and covering, it is not really
clear that we have settled on good measurements of ``complexity''.

\section{Complexity of sets}
\label{Complexity of sets}

	How might we measure complexity of sets?

	There are many answers to this question.  Here is
one.  Suppose that we have a set in some ${\bf R}^n$ which
is somehow ``$1$-dimensional''.  In terms of Hausdorff
measure for instance.  We can then ask whether it is contained
in a rectifiable curve, and if so how long the curve should be.

	The constraint that a set be contained in a rectifiable curve
imposes some nontrivial conditions on its size, but that is not the
only matter.  There is also a question of how {\em scattered} the
set is.

	Let us consider a concrete example, a Cantor set in the plane
constructed in the following manner.  On starts with the (closed) unit
square and then takes the four corner squares of one-fourth the size.
For each of these squares one does the same, and so forth.  At the
$k^{\rm th}$ stage of the construction one has $4^k$ squares of size
$4^{-k}$.  One takes the union of these $4^k$ squares at the $k^{\rm
th}$ stage to get a closed set $E_k$, and then one takes the
intersection of all the $E_k$'s to get a Cantor set $E$ which has
dimension $1$ in all the senses mentioned in the previous section.  In
particular, it satisfies (\ref{minkowski}) with $d = 1$, just like a
line segment.

	In terms of covering properties this set has practically
all the same properties as a line segment.  However one can show
that it is not contained in any curve of finite length, or even
in a countable union of them.  (See \cite{Mat}.)

	If a set is contained in a rectifiable curve, then one can use
the length of the shortest such curve as a measurement of its
complexity.  A measurement of ``infinite complexity'', because it
still allows infinitely many corners, loops, etc.  

	To be contained in a countable union of rectifiable curves
imposes already nontrivial restrictions on structure.  Roughly
speaking this property says that the set as approximate tangent lines
at almost all points.  See \cite{Fa, Fe, Mat} for details.  The Cantor
set described above is far from having approximate tangent lines.

	In \cite{J, O} there is a characterization given of the
compact sets $K$ which are contained in a curve of finite length, with
bounds for the length of the optimal curve.  This characterization is
given in terms of certain quantities which measure how close $K$ lies
to a line inside a given ball.  One looks at balls at all scales and
locations and combines them in a certain way to get a number which
turns out to be comparable to the length of the shortest curve which
contains $K$.  See \cite{J, O} for a more precise statement and further
details.  A key point is to take {\em orthogonality} into account, in
this case through the Pythagorean theorem.

	This story is related to the earlier one about functions.  We
can measure the complexity of a function through its oscillations,
like the ``modulus of continuity'' mentioned before.  Now it is more
natural to measure the complexity of a function in terms of the degree
to which it can be approximated by affine functions.  This is the
analogue of measuring the complexity of a $1$-dimensional set by the
extent to which it can be approximated by a line segment at different
scales and locations, and it is closely related to the topics in
Section \ref{Differentiability almost everywhere}.

	For either functions or sets one makes measurements in
arbitrary balls, say, and then combines all the measurements for all
balls to get some global measurement of complexity.  This is a very
basic theme in harmonic analysis, and it has natural geometric
counterparts as well.

	In this regard functions and sets are closely related.  If the
structure of a set is to be controlled through something like a
``parameterization'', as in the idea of putting a set inside a
rectifiable curve, then the parameterization involves functions whose
complexity as functions will be related to the geometric complexity of
the sets.  This will not necessarily work for all natural measurements
of the complexity of the functions, and not all of them will have
natural geometric meaning.

	See \cite{DS2, DS4, Secon, Seapp} for more information about
these themes, and also for analogous matters in higher dimensions.
For the ``qualitative'' versions (allowing countable unions of
surfaces) see \cite{Fa, Fe, Mat} also for higher dimensions.

	One should keep in mind that these are all measurements of
{\em infinite} complexity.  Even when a quantity like the length of a
curve is finite, {\em infinitely} many corners or loops or crossings
are permitted.  This is true of the other measurements that we have in
mind here as well.  The point is to control the ``level of infinity''.

	Of course these are not the only ways to think about
complexity of sets.  These are measurements of complexity in
comparison with {\em Euclidean geometry}, broadly construed.
Self-similar sets can also be considered as ``nice''.  They can admit
rapid computational description, for instance.  Sometimes
approximations of given sets by self-similar ones are even used for
the purposes of image compression.  See \cite{Ba}.

	The difference between roughly Euclidean geometry and models
which are self-similar but fractal is a bit like our earlier stories
for functions and operators.  One can make measurements of complexity
based on different methods of analysis -- the Fourier transform versus
wavelet bases, for instance -- and these can lead to very different
answers.

\section{Topological dimension}

	Let us now consider a different measurement of complexity of
sets, namely {\em topological dimension}.  \index{topological
dimension} The basic reference is
\cite{HW}.  For simplicity we shall restrict our attention to
compact metric spaces, but the theory works more generally.

	The definition of topological dimension in \cite{HW} is
inductive.  One starts by saying that a space has topological
dimension $0$ if it is totally disconnected.  A space has dimension
$\le n$ if each point has a system of neighborhoods whose boundary has
dimension $\le n-1$.  A space has dimension $n$ if it has dimension
$\le n$ but not $\le n-1$.

	I prefer another characterization.  A compact metric space
$(M, d(x,y))$ has dimension $\le n$ if for each $\epsilon > 0$ there
is a finite collection of open sets $\{U_i\}$ in $M$ such that the
$U_i$'s cover $M$, $\diam U_i < \epsilon$ for all $i$, and no point in
$M$ lies in more than $n+1$ of these sets.

	For example, an interval has dimension $\le 1$ but not
dimension $< 1$ unless it is a point.  This is not hard to check.  The
fact that a square has dimension $\le 2$ is slightly more amusing.
One has to be a bit careful in the way that one makes the covering,
but it is not difficult.

	It is important to notice what this characterization of topological
dimension does not ask for.  It does not ask for a bound on the total
number of the $U_i$'s.  In effect the condition controls the local
structure of the space but not at all how large it is.

	The condition also does not impose any requirements on the
shape of the $U_i$'s, only that they be sufficiently small.

	Of course the condition is genuinely {\em topological}, in the
sense that it does not depend on the choice of metric, only on the
underlying topology.  Change the metric and the space might be much
``larger'' in terms of size, but not in terms of the topological
dimension.  For instance, if the original metric is $d(x,y)$, and we
change it to $d(x,y)^\frac{1}{2}$, then this is still a metric -- the
main point is that the triangle inequality continues to hold -- and
dimensions which measure size (like the Hausdorff dimension) increase
by a factor of $2$ under this change, but the topological dimension
remains the same.  Similar results hold if one replaces $\frac{1}{2}$
with any number $s$ so long as $0 < s < 1$.

	There is a nice theorem of Alexandroff, \index{Alexandroff's theorem}
to the effect that a
compact metric space $M$ has topological dimension $\le n$ if and only
if it can be approximated by finite polyhedra of dimension $\le n$ in
the following sense.  Given an arbitrary $\epsilon > 0$ one asks for a
a finite polyhedron $P$ and a continuous mapping $f : M \to P$ with
the property that each ``fiber'' $f^{-1}(p)$, $p
\in P$, has diameter $< \epsilon$.  This does not imply that $M$ is
homeomorphic to a subset of finite polyhedron of dimension $\le n$,
but it wants to approximate that statement.  For instance, suppose
that $M$ is a subset of a plane, the Cartesian product of the unit
interval $[0,1]$ and a totally disconnected infinite compact set, like a
convergent sequence together with its limit.  It is easy to see that
this set has topological dimension $\le 1$, but of course it is not
homeomorphic to a subset of a $1$-dimensional polyhedron.

	Note well that the polyhedron $P$ is permitted to depend on
$\epsilon$.  One can also look at its complexity, total number of
faces, for instance.  The notion of topological dimension does not
impose any restriction on the complexity of $P$, except for its
dimension, but one can profitably impose such conditions.  See
\cite{Setop}, for instance.

	Another nice theorem says that if $M$ has topological
dimension $\le n$ then $M$ can be homeomorphically embedded
into ${\bf R}^{2n+1}$.  Thus the notion of topological dimension
provides a characterization of the compact metric spaces which
can be embedded into some finite-dimensional Euclidean space.

	This embedding result is quite remarkable, in the way that it
involves infinite processes.  The direct part of Alexandroff's
theorem, the passage from the bound on the topological dimension to
the polyhedral approximation, is a simpler matter of finite
combinatorics and partitions of unity.  One fixes an $\epsilon > 0$,
applies the hypothesis at a particular scale, obtains the conclusion
at a particular scale.  By contrast the embedding theorem involves
distances at all scales; one has to build a mapping for which the
inverse image of a point is a point, and not merely small in diameter.
The usual proof relies on the Baire category theorem.

	The embedding theorem provides a nice example of a result
which involves infinite processes and for which it is convenient
not to make the proof at all explicit!

	Another way to understand topological dimension involves the
notion of {\em stable values}.  \index{stable values} Let $M$ be a
compact metric space, and let $I_n$ denote the closed unit cube in
${\bf R}^n$.  We say that a point $p$ in $I_n$ is a {\em stable
value} of a continuous mapping $f : M \to I_n$ if there is a $\delta >
0$ so that if $g : M \to I_n$ is any continuous mapping which is
$\delta$-close to $f$, i.e.,
\begin{equation}
	|f(x) - g(x)| < \delta \qquad\hbox{for all } x \in M,
\end{equation}
then $p$ lies in the image of $g$.  Of course this implies that
$p$ lies in the image of $f$, and we call $p \in f(M)$ an
{\em unstable value} of $f$ if it does not enjoy this property.

	It turns out that if $M$ has topological dimension strictly
less than $n$ then every value of every continuous mapping
$f : M \to I_n$ is unstable, and that if $M$ has topological dimension
at least $n$ then there is a continuous mapping $f : M \to I_n$
with a stable value.  (One can take the stable value to be the
center of the cube $I_n$.)

	In other words, if the topological dimension is strictly less
than $n$, then $M$ is too ``thin'' to have a mapping into $I_n$ which
is very ``serious''.  Think about Peano curves, which can provide
examples of continuous mappings from a line segment onto a square.  Or
the Cantor function, which maps the Cantor set continuously onto an
interval.  It is easy to map spaces of lower dimension onto spaces of
larger dimension, but these mappings can be approximated by mappings
with much smaller image.  This is very easy to see for the examples
just mentioned, and indeed they are typically obtained as the limit
of such mappings.

	The fact that spaces of dimension at least $n$ admit mappings
to cubes of dimension $n$ with stable values is very nice.  It is a
strong way of saying that the space really has some $n$-dimensional
flesh.  This result requires infinite processes, and the usual proof
uses the Baire category theorem.

	There are other results to this effect, that having
topological dimension $< n$ means that the space is approximately
``invisible'' for topological considerations of dimension $n$ or
greater, in terms of homotopy properties of mappings into spheres, and
in terms of homology and cohomology.  See \cite{HW} for details.

	One might think of the {\em topological dimension} as being
analogous to the notion of {\em depth} as a measurement in complexity
theory.  Let us think about making constructions on a fixed polyhedron
$P$ of dimension $n$, constructions of a continuous mapping with some
properties, say.  A common line of reasoning begins by making choices
at the vertices, then the edges, then the $2$-dimensional faces, and
so forth through dimensions until one reaches the top dimension $n$.
It happens sometimes that the constructions on the different faces of
a given dimension are more-or-less independent of each other, that the
more tricky point is the number of layers, the ``depth'', which in
this case means the dimension.  Thus the complexity of the polyhedron
-- the total number of constituent simplices for instance -- is
sometimes not so crucial, while the dimension is.  Although this is
simpler for polyhedra than general spaces, similar considerations
apply in general, with Alexandroff's approximation theorem providing a
convenient bridge.

	This type of phenomenon comes across nicely in \cite{P2}.
See also \cite{Setop}.

	The Brouwer whose name appears so frequently and prominently
in \cite{HW} is the same as the eminent logician.

\section{Topological dimension, size, and structure}

	The topological dimension is a very good measurement of
``size'' of a space for looking at the behavior of continuous mappings
and other topological questions.  In other ways it provides
practically no information; metrically the space could be very large.
We can make this precise with the concept of Hausdorff measure and
Hausdorff dimension.  We defined these in Section
\ref{Sizes of sets} only for subsets of Euclidean spaces, but the
definitions extend to metric spaces in a straightforward way.

	A metric space can have topological dimension $0$ but large
Hausdorff dimension, or even infinite Hausdorff dimension.  One can
make versions of the usual Cantor set construction within spaces of
larger dimension, or even infinite-dimensional spaces, to get
examples.  One can also proceed more abstractly.  Suppose that $(M,
d(x,y))$ is a metric space of Hausdorff dimension $d$.  We can apply
the {\em snowflake transform} to it to get the space $(M, d(x,y)^s)$,
where $0 < s < 1$.  This is again a metric space (the triangle
inequality remains valid), and it has Hausdorff dimension $d/s$.

	Of course the snowflake transform does not change the
underlying topological structure, or the topological dimension.  It is
irrelevant for purely topological considerations such as the behavior
of continuous mappings, the existence of stable values, etc.  And yet
it affects the {\em geometry} of the space in a dramatic way.  In
terms of how many balls of a given radius are needed to cover the
space, for instance.

	In this regard topological dimension is rather remarkable
for extracting a {\em topological} measurement of size which does
not worry too much about {\em geometric} measurement of size.

	Is there a relationship between these measurements of size?
Is the geometric one always larger than the topological one, for instance?

	The answer is yes.  The topological dimension is always less than
or equal to the Hausdorff dimension.  Note well that the topological
dimension is necessarily an integer, while the Hausdorff dimension
need not be.

	This result can be found in \cite{HW}.  In fact there is a
more precise statement: if $M$ has topological dimension at least $n$,
then $M$ has {\em positive} Hausdorff measure in that dimension as
well.  (These results in \cite{HW} for Hausdorff dimension and
positive Hausdorff measure are actually stated in terms of Hausdorff
{\em content}, but this is equivalent, as in the remarks of Section
\ref{Sizes of sets}.)

	What about a converse?  Topological dimension need not control
Hausdorff dimension for a {\em particular} metric, but topological
dimension is independent of the choice of metric.  Can one always find
a metric which is compatible with the topology and for which the
Hausdorff dimension is equal to the topological dimension?  The answer
is yes, and one can even embed a given compact metric space of finite
topological dimension into a finite-dimensional Euclidean space in
such a way that the image has Hausdorff dimension equal to the
topological dimension, where Hausdorff dimension is now computed using
the ambient Euclidean metric.  See \cite{HW}.

	It may not be possible to find a metric on a given space so that
the space has {\em finite} Hausdorff measure in the topological dimension.
One can make easy counterexamples using the Cartesian product of an
interval with an infinite set.

	In summary, topological dimension and Hausdorff dimension both
provide measurements of the structure of a set.  They are both
measurements of ``infinite'' levels of complexity; they tolerate
infinite complexity but control it in some ways and not in other ways.
They provide very different measurements of complexity.

	The ideas of Section \ref{Complexity of sets} point to another
measurement of complexity still.  A measurement that looks for a {\em
geometric} resemblance to Euclidean spaces rather than a topological
one, or mere similarity in terms of rough measurements of size.

	Hausdorff dimension is always at least as big as topological
dimension, what happens when they are equal?  The topological
dimension provides a kind of ``floor'' for the geometric size, but if
the actual size is not much larger, then can we expect some nontrivial
restrictions on the geometric structure?  Results which permit one to
obtain stronger conclusions about the geometric complexity of the
space?  There are some results of this nature, and this general
question is discussed in \cite{Secon}.

\section{Lipschitz mappings}
\label{Lipschitz mappings}

	A {\em Lipschitz mapping} between metric spaces is one which
does not expand distances by more than a constant factor.

	More formally, let $(M, d(x,y))$ and $(N, \rho(u,v))$ denote
two metric spaces.  We say that a mapping $f : M \to N$ is {\em Lipschitz}
if there is a constant $C > 0$ so that
\begin{equation}
	\rho(f(x), f(y)) \le C \, d(x,y) \qquad\hbox{for all } x,y \in M.
\end{equation}
We might say {$C$-Lipschitz} to make the constant more precise.

	It will be convenient to use the phrase {\em Lipschitz function}
to mean a Lipschitz mapping into the real line (with the standard metric).

	The Lipschitz condition provides a measurement of infinite
complexity of functions.  It allows infinite amounts of activity, even
for functions on the unit interval, but nonetheless imposes nontrivial
restriction on this activity.

	To put this measurement into perspective it is helpful to
compare the Lipschitz condition with our earlier discussion of
complexity of functions, and to view our earlier discussion with new
eyes of suspicion.  The Lipschitz property has the nice feature that
it is very geometric and simple, making sense for mappings between
arbitrary metric spaces.  Most of the earlier discussion was not so
robust, requiring much simpler structure for both domain and range.
Think of approximations by polynomials, rational functions, splines,
trigonometric polynomials, etc.  All of these require more precise
structure of the underlying space even to be formulated.

	Continuous (real-valued) functions can always be approximated
by Lipschitz functions, modulo some technical hypotheses.  Let $(M,
d(x,y))$ be a metric space, and suppose for a moment that $M$ is
compact.  Let $f : M \to {\bf R}$ be a continuous function, and set
\begin{equation}
	f_j(x) = \inf \{f(y) + j \, d(x,y) : y \in M\}.
\end{equation}
This defines a real-valued function on $M$, and one can show that
$f_j$ is $j$-Lipschitz.  This is not hard, using the fact that $x
\mapsto d(x,y)$ is $1$-Lipschitz for each $y \in M$, which itself
comes from the triangle inequality.  One can also prove that $f_j$
converges to $f$ uniformly on $M$ as $j \to \infty$.  This is well
known and not very difficult to verify, and we omit the details.  The
main point is to use the fact that $f$ is uniformly continuous since
it is continuous and $M$ is compact.  Indeed one can relate the rate
of convergence of $f_j \to f$ to the modulus of continuity of $f$ in a
simple way.

	This argument works on general metric spaces if one makes some
additional assumptions.  If $M$ is not compact then the finiteness of
$f_j$ is not assured, but it can be verified under moderate
assumptions on $f$ (i.e., its growth at infinity).  Once one knows
that $f_j$ is finite at any point one can show that it is finite
everywhere and $j$-Lipschitz.  Uniform convergence to $f$ can be proved
under the assumptions of uniform continuity of $f$ and restrictions
on the growth of $f$ at infinity.  (Boundedness would be plenty.)

	Under modest conditions on the complexity of $M$ one can make
other procedures of approximation which are linear in $f$ and more
obvious ``locality'' properties.  This is discussed in
\cite{Seapp}.

	Lipschitz mappings cooperate well with our measurements of
{\em sizes} of sets.  Let $(M, d(x,y))$ and $(N, \rho(u,v))$ be two
metric spaces, and suppose that $f : M \to N$ is $C$-Lipschitz.
Let $E$ be any subset of $M$.  Then the size of $E$ can be related
to the size of $f(E)$ in the following ways.  Given a dimension $d \ge 0$,
we have that
\begin{equation}
	H^d(f(E)) \le C^d \, H^d(E) \quad\hbox{ and }\quad
			A^d(f(E)) \le C^d \, A^d(E).
\end{equation}
This is easy to check from the definitions.  Similarly, if
$N(E, \epsilon)$ is defined as in the beginning of Section \ref{Sizes of
sets}, as the smallest number of balls of radius $\epsilon$ needed
to cover $E$, then we have that
\begin{equation}
	N(f(E), C \, \epsilon) \le N(E, \epsilon).
\end{equation}

	A particular consequence of these inequalities is that
Lipschitz mappings never increase Hausdorff dimension.  However
they can increase topological dimension.  To see this, let us first
define a ``symbolic'' version of the Cantor set.  We take $K$
to denote the set of all infinite binary sequences, i.e., all
sequences $\{x_j\}_{j=1}^\infty$ such that $x_j \in \{0,1\}$
for all $j$.  We equip this set with a metric by setting
$D(\{x_j\}, \{y_j\}) = 0$ when $x_j = y_j$ for all $j$, 
\begin{equation}
	D(\{x_j\}, \{y_j\}) = 2^{-l} \quad\hbox{when } x_j = y_j
			\hbox{ for all } j \le l, x_{l+1} \ne y_{l+1}.
\end{equation}
It is easy to see that this defines a metric on $K$ in such a way that
$K$ is compact.  One can also check that $K$ has Hausdorff dimension $1$.

	This metric is substantially different from the one on the
standard middle-thirds Cantor set which is induced by the usual
Euclidean metric on the line, because they give different Hausdorff
dimensions.  The standard metric for the middle-thirds Cantor set
corresponds approximately to a power of this metric.

	Define a mapping $f : K \to [0,1]$ by
\begin{equation}
	f(\{x_j\}) = \sum_{j=1}^\infty x_j \, 2^{-j}.
\end{equation}
One can check that this mapping is $1$-Lipschitz and also surjective.
However, $K$ has topological dimension $0$, while $[0,1]$ has
topological dimension $1$.  Thus Lipschitz mappings can indeed
increase the topological dimension.

	Of course Lipschitz mappings can decrease topological
dimension, e.g., by mapping $[0,1]$ to a point.  One cannot map
$[0,1]$ into a set like $K$ in an interesting way, the image will
always have to be a point.  In situations like this -- continuous
mappings between compact spaces -- in order for the topological
dimension of the image to be less than that of the domain it is
necessary that there be a point in the range whose inverse image has
topological dimension at least the difference between the topological
dimensions of the domain and image.  See Theorem VI 7 on p91 of
\cite{HW}.

\section{Differentiability, rectifiability, etc.}

	There are issues of infinite complexity for Lipschitz mappings
which are much more subtle than mere size.  A Lipschitz function on a
Euclidean space is differentiable almost everywhere.  Thus one does
not simply have bounds on the sizes of the oscillations of Lipschitz
functions, but also a kind of asymptotic rigidity at almost all
points.  See \cite{Fe, St1, Seapp}.  A well-known variant of this
famous result says that given any $\epsilon > 0$ one can find a $C^1$
function which agrees with the given Lipschitz function except on a
set of measure $< \epsilon$.  (See \cite{Fe}.)

	There are also quantitative versions of these rigidity
properties of Lipschitz functions on Euclidean spaces, as discussed in
\cite{Do, DS4, Secon, Seapp}.  (See also Section \ref{Big pieces of 
bilipschitz mappings} below.)

	To put these matters into perspective, consider the {\em
H\"older continuous} functions on ${\bf R}^n$.  That is, $f : {\bf
R}^n \to {\bf R}$ is {\em H\"older continuous of order $\alpha$} if
there is a constant $C > 0$ so that
\begin{equation}
	|f(x) - f(y)| \le C \, |x-y|^\alpha
\end{equation}
for all $x, y \in {\bf R}^n$.  This is the same as the Lipschitz
condition when $\alpha = 1$.  When $\alpha > 1$ it is too strong, and
implies that $f$ must be constant.  When $\alpha < 1$ it is weaker
than the Lipschitz condition in the way that it controls the local
behavior of $f$.  One can see that there is nothing like the above
rigidity properties of Lipschitz functions for H\"older continuous
functions of order $\alpha < 1$.  This is discussed in more detail in
\cite{Seapp}.

	Similarly, for general metric spaces one cannot hope to have
the kind of rigidity properties for Lipschitz functions that one has
on Euclidean spaces.  Lipschitz functions on our symbolic Cantor set
$K$ of Section \ref{Lipschitz mappings}, for instance, do not enjoy
such rigidity properties.  Nor do Lipschitz functions on Euclidean
spaces equipped with ``snowflake'' metrics of the form $|x-y|^s$, $0 <
s < 1$.  (A Lipschitz function on such a snowflake space is
essentially the same thing as a H\"older continuous function with
respect to the ordinary Euclidean metric.)  There are some metric
spaces which are far from Euclidean in their geometry but for which
one still has a version of differentiability almost everywhere, namely
Heisenberg groups with their Carnot metrics.  See \cite{Pa}.

	Thus Lipschitz functions are very interesting in terms of
infinite complexity.  They cooperate with measurements of size in an
obvious way, but they can also enjoy additional rigidity properties
which are not explicitly part of their definition.  These rigidity
properties do not work for metric spaces in general, but depend on the
geometry of the underlying space.  These ``hidden'' rigidity properties
are absent also for many other conditions that one might impose
on functions, like H\"older continuity.

	Just as one had differentiability (almost everywhere) and
other natural ``rigidity'' properties of functions, there are
geometric conditions like {\em rectifiability} and {\em uniform
rectifiability} for sets.  In the geometric setting one loses much of
the linear structure of functions and methods for analyzing them (like
the Fourier transform and wavelets), but one gains a more rich
language, so that there are completely different types of criteria for
the presence of good structure.  Criteria involving the relationship
between size and topology, for instance.  See \cite{Fa, Fe, Mat, DS4,
Secon, Seapp}.

	For further discussion about the analogy between functions and
sets and the methods for analyzing them see \cite{DS3, DS4}, and note
the differences in their conclusions!

	See \cite{DS6} for a search for other types of rigidity
properties related to Lipschitz mappings and geometric measure theory.

	The preceding paragraphs are very sketchy, little more than
pointers to a sizable literature.  In short the main point is that in
these directions there are worlds of situations in which infinite
complexity arises naturally and nontrivial ways to measure it and to
distinguish between different levels of complexity.

\section{Wild embeddings}
\label{Wild embeddings}
\index{wild embeddings}

	Given an arc in a Euclidean space, can we always straighten it
by a global homeomorphism?

	By an {\em arc} we mean a set which is homeomorphic to the
unit interval $[0,1]$.  We might say ``topological arc'' for emphasis.
The question asks whether we can find a global homeomorphism from the
ambient Euclidean space onto itself which sends the arc to a straight
line segment.

	In ${\bf R}^1$ this is trivially true because all arcs are
line segments.  In ${\bf R}^2$ this is true but nontrivial, one of the
standard facts about plane topology.  Starting in ${\bf R}^3$ it fails
to be true.  See \cite{Moi}, for instance.

	Similarly, given a compact set in ${\bf R}^n$ which is
homeomorphic to the usual Cantor set, one can ask whether there is a
global homeomorphism on ${\bf R}^n$ which takes the given set to a
subset of a line.  Again the answer is yes in dimension $2$ and no
starting in dimension $3$.  Basic examples are provided by ``Antoine's
necklaces'', which are topological Cantor sets in ${\bf R}^3$ with
non-simply connected complement.  If there were a global homeomorphism
carrying such a set to a subset of a line, then the complement would
have to be simply-connected (because the set is totally disconnected).

	The construction of the Antoine's necklace is easy enough to
describe.  For the usual Cantor set one takes a closed interval,
replaces it by two disjoint closed subintervals, performs a similar
operation on those, etc.  For the Antoine's necklace one starts with a
solid torus and replaces it by a chain of solid tori contained inside.
These new solid tori are disjoint but linked.  In each of these solid
tori one repeats the construction, and does this forever.  In the
limit one gets a set which is homeomorphic to a Cantor set, but the
complement is not simply connected, as one can verify.  See
\cite{Moi} for details.

	Similarly there are homeomorphic embeddings of the $2$-sphere
into ${\bf R}^3$ which are ``wild'', i.e., cannot be transformed into
a standard $2$-sphere by a global homeomorphism.  See \cite{Moi}
again, and also \cite{BuC, Dv1}.

	In all of these situations one has a subset of a Euclidean
space which is homeomorphic to something fairly simple but for which
the embedding has to have an infinite amount of twisting.  These
phenomena also lend themselves naturally to questions of complexity.
One can for instance look at finite versions of the constructions and
measure the amount of distortion needed for global homeomorphisms
which make the straightenings.  The impossibility of straightening the
infinite constructions implies the absence of uniform bounds for the
finite versions.  (Uniform bounds on the moduli of continuity of the
relevant homeomorphisms and their inverses, that is.)

\section{The Whitehead continuum}
\index{Whitehead continuum}

	Is every contractible manifold homeomorphic to a ball?  

	In dimensions $1$ and $2$ this is true and classical.  It is
not true starting in dimension $3$.  The first counterexample was
obtained by taking the complement of the {\em Whitehead continuum} in
the $3$-sphere.  Roughly speaking, the {\em Whitehead continuum} is
defined as follows.  One starts with a (compact) solid torus $T_1$,
which we think of as embedded in the $3$-sphere.  That is our first
approximation.  For the second approximation we take another solid
torus $T_2$ which is embedded inside $T_1$.  The embedding is rather
special; $T_2$ is stretched around the middle of $T_1$ and has its
ends clasped together on the other side.  (See \cite{Dv2, K} for
pictures.)  This clasping makes it impossible to shrink $T_2$ to a
point without leaving $T_1$ and without having a copy of $T_2$ cross
itself in the course of the shrinking.  On the other hand $T_2$ can be
shrunk to a point inside $T_1$ if one does allow copies of $T_2$ to
cross themselves in the shrinking (but still not to ever leave $T_1$).
This is because of the way that $T_1$ is stretched around the center of
$T_2$; it does not simply go around in a circle, and in particular
does not contain a homotopically nontrivial loop in $T_1$.

	$T_2$ is only the second approximation to the Whitehead
continuum.  For the third we embed a new solid torus inside $T_2$ in
the same way that $T_2$ was embedded into $T_1$.  This procedure is
repeated indefinitely, at each stage a new solid torus is embedded
into the last one in the same way that $T_2$ was embedded into $T_1$.
To get the Whitehead continuum one takes the intersection of this
decreasing chain of embedded tori.  Let us call the final result $W$.
This is a compact set, and it is easy to see that it is connected.

	In fact, ${\bf S}^3 \backslash W$ is {\em contractible}.  This
means that it can be continuously deformed inside itself until it
becomes a single point.

	${\bf S}^3 \backslash W$ is also an open subset of ${\bf
S}^3$, hence a (noncompact) manifold without boundary.  It is
contractible but not homeomorphic to a $3$-dimensional open ball.  The
reason is that it is not {\em simply-connected at infinity}.  To
explain what this means it is helpful to call an open subset of ${\bf
S}^3 \backslash W$ a {\em neighborhood of infinity} if it is the
complement in ${\bf S}^3 \backslash W$ of a compact subset of ${\bf
S}^3 \backslash W$.  If ${\bf S}^3 \backslash W$ were simply-connected
at infinity it would mean that for each neighborhood $U$ of infinity
there is a smaller neighborhood $V$ such that every loop in $V$ can be
contracted to a point inside $U$.  This is no true for ${\bf S}^3
\backslash W$.  See \cite{Dv2, K}.

	Similar issues arise in the wild examples mentioned in the
previous section.  Sometimes the problem is more blatant, e.g., a
region is not even simply connected when it ``should'' be, like the
complement of Antoine's necklaces.  In all cases there is an
``asymptotic'' problem, a problem with homotoping loops in the
complement of the set but which lie near the set.

\section{Decompositions}
\index{decomposition spaces}

	Imagine taking ${\bf R}^3$ and some compact set $K \subseteq
{\bf R}^3$ and shrinking $K$ to a point while leaving the other points
alone.  What happens to the topology?

	If $K$ is a straight line segment the resulting space is {\em
homeomorphic} to ${\bf R}^3$.  The operation of collapsing the line
segment to a point is pretty severe geometrically, but it does give
back the same topological space.

	This is not true if $K$ is a circle.  If we collapse a circle
to a point then the complement of that point in the resulting space is
the same topologically as the complement of the circle in ${\bf R}^3$,
and is not simply-connected.  Since the complement of a point in ${\bf
R}^3$ is always simply-connected, the new space that we made cannot be
homeomorphic to ${\bf R}^3$.

	If we take $K$ to be the Whitehead continuum (taking care to
have it inside ${\bf R}^3$ and not just ${\bf S}^3$), then the
resulting space is again different from ${\bf R}^3$ topologically.
Indeed, every point in ${\bf R}^3$ has arbitrarily small punctured
neighborhoods which are simply connected, and this is not true when we
collapse the Whitehead continuum to a point.  This corresponds to the
failure of simple-connectivity at infinity described in the preceding
section.

	Similarly if one collapses a wild arc to a point one may not get
${\bf R}^3$ back again.

	In these cases one not only does not get ${\bf R}^3$ back
again, one does not even get a topological manifold, because of the
behavior near the strange point.

	More generally one can collapse families of compact sets to
points rather than just individual ones.  One can simply speak of {\em
decompositions} of ${\bf R}^3$ (or other spaces), in which one takes a
partition (by closed sets say) and uses this to make a quotient
topological space.  This is discussed in more detail in \cite{Dv2}.
There are natural conditions to impose to avoid pathologies as one
gets by collapsing a circle or a Whitehead continuum to a point, but
there is not a simple way to tell when one will get back ${\bf R}^3$
or something strange.

	What about natural finite approximations to making spaces by
collapsing subsets to points?  Can one make measurements of
complexity?  One way to do this is given in \cite{Seqs}.  There one
considers decompositions obtained through a certain kind of iteration,
for which the Whitehead continuum is a simple special case.  The
compact sets to be collapsed are obtained as limits of compact domains
obtained by repeated embeddings (roughly speaking).  This kind of
construction has a natural {\em topological} ``self-similarity'', to
which one associates a natural geometry which is self-similar in the
usual sense.  In this geometry the compact sets to be collapsed to
points have diameter zero, as they should, but their approximations
are shrinking at a particular rate.  The final spaces can be viewed as
subsets of ${\bf R}^4$, with the induced geometry.

	In making these geometries it makes sense to stop the
construction after finitely many generations.  These approximations
are automatically homeomorphic to ${\bf R}^3$, the question is what
kinds of estimates one can have for the homeomorphisms, in terms of
the approximately self-similar geometry.  Some observations and
examples are given in \cite{Seqs}, but these quantitative issues have
yet to be studied much.

\section{Domains in Euclidean spaces}
\label{Domains in Euclidean spaces}

	Let $\Omega$ be an open subset of some ${\bf R}^n$.  Let us
define the {\em quasihyperbolic metric} \index{quasihyperbolic metric}
on $\Omega$ to be the Riemannian metric given at an arbitrary point $x
\in \Omega$ by the formula
\begin{equation}
\label{qh metric}
	\frac{ds^2}{\dist(x,\partial \Omega)^2}.
\end{equation}
Thus the length of an infinitesimal bit of curve at $x$ is the same as
the ordinary Euclidean length divided by the distance to the boundary
of $\Omega$.

	This defines a metric which blows up as one goes to the boundary.
Not only does the infinitesimal form blow up, but the length of any
path which goes toward $\partial \Omega$ will blow up.  The approximate
geometry is captured well by the following observation.  In the geodesic
distance for the quasihyperbolic metric, a ball with center $x \in \Omega$
and radius $1$ is approximately like a Euclidean ball with center $x$
whose radius is proportional to $\dist(x,\partial \Omega)$ by a positive
constant that stays away from $1$ by a definite amount.

	One might just as well define the quasihyperbolic metric using
a smooth function $\delta_\Omega(x)$ instead of $\dist(x,\partial
\Omega)$, but where $\delta_\Omega(x)$ is chosen so that it is bounded
from above and below by positive multiples of $\dist(x,\partial
\Omega)$ (as on Theorem 2 on p171 of \cite{St1}).  

	If the complement of $\Omega$ is compact one might prefer to
adjust this definition, working on the sphere instead so that the
point at infinity is included, etc.  Of course one can extend this
construction to domains on smooth manifolds.

	Using this metric one can look at the (Euclidean) boundary
behavior of $\Omega$ as asymptotic geometry at infinity with respect
to the quasihyperbolic geometry.  If $\Omega$ is bounded by a compact
smooth hypersurface, say, then the quasihyperbolic metric looks a lot
like the usual models for hyperbolic geometry.  Otherwise it can be
much more complicated.

	This geometry is often natural for analysis.  Bounded harmonic
functions on $\Omega$, for instance, have very tame behavior on
quasihyperbolic balls of bounded radius, while their behavior near the
boundary can be much more complicated.

	If the boundary of $\Omega$ is itself complicated, the
quasihyperholic metric provides a natural way to think about exploring
it.  Imagine that $\Omega$ is a sea, and that were are in a small boat
exploring the fjord which borders it.  In the examples of Section
\ref{Wild embeddings} there is interesting behavior already on
bounded subsets of $\Omega$ with respect to the quasihyperbolic
metric.  There can be loops which cannot be filled by disks within a
bounded region (in the quasihyperbolic metric), and this behavior can
be repeated infinitely often as one approaches the boundary of
$\Omega$.  One way to make interesting measurements is to count the
number of times that such events occur, e.g., how numerous they are
within regions of a certain size.  This idea works well in the context
of ``uniformly rectifiable sets'' and ``almost flat hypersurfaces''
(Section \ref{Almost flat hypersurfaces}), for which there are
typically infinitely many ``complications'' whose frequency can be
counted with fairly good precision.  See \cite{DS2, DS4, Seag, casscI,
casscII, pams, Seapp}.

	Alternatively one can look at all locations in $\Omega$ in a uniform
way but try to deal in a more subtle way with the range and nature of the
complications.  See \cite{Al, AV1, AV2, HY, V}.

\section{Spaces of bounded geometry}

	Let $(M, g)$ be a complete Riemannian manifold of dimension
$n$.  Let us say that $(M, g)$ has {\em bounded geometry} if for each
$x \in M$ the ball in $M$ with center admits a diffeomorphic mapping
onto an open subset of ${\bf R}^n$ in such a way that this
diffeomorphism distorts distances by only a bounded factor and its
derivatives are uniformly bounded.  We do not allow these bounds to
depend on the point $x$, but the bounds on the derivatives are
permitted to depend on the level of differentiation.

	In a space of bounded geometry the local behavior is
controlled in a strong way, and the issue is what happens at infinity.

\beginexamples
{\rm 

	(1) Euclidean and hyperbolic spaces have bounded geometry. 

	(2)  Homogeneous spaces always have bounded geometry.

	(3) Domains in Euclidean spaces equipped with their
quasihyperbolic metrics (as in Section \ref{Domains in Euclidean
spaces}) are spaces of bounded geometry, at least if one takes the
trouble to replace $\dist(x,\partial \Omega)$ in (\ref{qh metric}) by
a smooth version of the distance, as discussed above.

	(4) The universal covering of any compact Riemannian manifold
(without boundary) has bounded geometry.  

	(5) A leaf of a smooth foliation in a compact manifold has
bounded geometry (through the induced metric).
}
\end{examples}

	These examples are in a sense very different from each other,
and they have lives of their own, but they can be profitably seen to
some extent as special cases of a general notion.  They also make for
good examples of infinite complexity within ordinary mathematics.
There are many natural ways to measure this asymptotic complexity,
including isoperimetric inequalities.  Some references include
\cite{Gr0, Gr2, Gr3, At1, At2, AtH}.

\chapter{Some results about geometric complexity}

	In the preceding chapters we saw a range of situations in
which one can make measurements of complexity even when the complexity
is infinite.  In this chapter we discuss a few particular results of
this nature.

\section{Discrete parameterizations}
\label{Discrete parameterizations}

	The following result is proved in \cite{Sto}.

\begintheorem 
\label{stong's theorem}
Let $r$ and $s$ be positive integers, with $s \ge r$.  Then there is a
{\em bijection} $f : {\bf Z}^r \to {\bf Z}^s$ such that
\begin{equation}
\label{holder bound for the mapping}
	|f(x) - f(y)| \le C \, |x-y|^\frac{r}{s}
\end{equation}
for some constant $C > 0$ and all $x, y \in {\bf Z}^r$.
\end{theorem}

	This is quite nice.  We are all familiar with space-filling
curves, we know that they have to be somehow ``fractal'', but that is
not at all the same as getting the discrete version to work out right,
on the nose.

	The H\"older exponent $\frac{r}{s}$ in (\ref{holder bound for
the mapping}) is absolutely the right one.  It exactly matches the
differences in the ``dimensions'' of ${\bf Z}^r$ and ${\bf Z}^s$.
Notice, for instance, that the number of elements of ${\bf Z}^r$ which
lie in a ball of radius $L > 1$ centered at the origin is bounded from
above and below by constant multiples of $L^r$, and there is a similar
statement for ${\bf Z}^s$.  One can derive the sharpness of the
H\"older exponent $\frac{r}{s}$ from this observation.

	To put the theorem into perspective it is helpful to think
more about its continuous analogues, mappings from ${\bf R}^r$ to
${\bf R}^s$.  For this one cannot ever get bijectivity when $s > r$,
by standard results in topology.  Bijectivity can only occur in a
discrete setting, but then what does it mean exactly?

	One can think of bijectivity as a way to say
``measure-preserving'', for instance.  It implies the following
``regularity property'', which arose first in \cite{D3}.
\index{regular mappings}
If $f$ is as in the theorem, then there is a constant $K > 0$ so that
for each ball $B$ of radius $t$ in ${\bf R}^s$ we can cover the
inverse-image $f^{-1}(B)$ by $\le K$ balls of radius $t^\frac{s}{r}$.
This observation is analogous to one in
\cite{D3}, and can be proved as follows.  Let $t$ and
$B$ be given as above. We may as well assume that $t \ge 1$.  Let $A$
be any subset of $f^{-1}(B)$ with the property that $|x-y| \ge
t^\frac{s}{r}$ for all $x, y \in A$ with $x \ne y$.  We want to derive
a uniform bound for the possible number of elements in $A$.  Given $x
\in A$, let $B(x)$ denote the ball $B(x, \frac{1}{2} \, t^\frac{s}{r})$.  
The balls $B(x)$, $x \in A$, are pairwise disjoint, and so their
images are pairwise disjoint as well, since $f$ is injective.  On the
other hand their images $f(B(x))$, $x \in A$, are all contained in the
ball obtained by dilating $B$ by a constant factor (which depends on
$f$ but not on $B$).  This is easy to check, since we know that $f(x)$
lies in $B$ for each $x \in A$, and since the H\"older condition on
$f$ implies that each $f(B(x))$ has diameter bounded by a constant
multiple of $t$.  We know that each ball $B(x)$ contains about
$(t^\frac{s}{r})^r = t^s$ elements of ${\bf Z}^r$ (to within a
constant factor), and so the number of elements of $\bigcup_{x \in A}
f(B(x))$ is at least a constant multiple of $t^s$ times the number of
elements of $A$.  All these points must lie within a fixed dilate of
$B$, and this implies that $A$ contains no more than a bounded number
of elements, with a bound that does not depend on $B$.  This implies
the bounded covering property that we want.
	
	The measure-preserving and regularity properties make sense
for mappings between ${\bf R}^r$ and ${\bf R}^s$ as well.  In fact,
the ``discrete'' result of Theorem \ref{stong's theorem} implies
continuous versions through normal families arguments, as follows.
Fix $r$ and $s$ as above, and choose $f$ as in Theorem \ref{stong's
theorem} and so that $f(0) = 0$.  We can always arrange this
normalization.  Given $\epsilon > 0$, let $f_\epsilon$ be defined by
\begin{equation}
	f_\epsilon(x) = \epsilon^\frac{r}{s} \, f(\frac{x}{\epsilon}),
\end{equation}
so that $f$ defines a bijection between $\epsilon \cdot {\bf Z}^r$ and
$\epsilon^\frac{r}{s} \cdot {\bf Z}^s$.  Each $f_\epsilon$ satisfies
the analogue (\ref{holder bound for the mapping}), with the same
constant, because we were careful about the scaling.  The
Arzela-Ascoli argument permits one to find a sequence of $\epsilon$'s
tending to $0$ such that the corresponding $f_\epsilon$'s ``converge
uniformly on compact subsets'' to a mapping $g : {\bf R}^r \to {\bf
R}^s$.  There is a small technical problem here, since the
$f_\epsilon$'s do not all have the domain, but this is not a serious
problem, because the domains become increasingly thick.  This kind of
issue is discussed in some detail in \cite{DS6}.  The mapping $g$ that
results might be considered as a continuous version of the mapping $f$
from Theorem \ref{stong's theorem}, since it can be derived from $f$
in this way.

	Our limiting mapping $g$ also satisfies the analogue of
(\ref{holder bound for the mapping}), since the $f_\epsilon$'s do, and
with a uniform bound.  One can also check that $g$ is {\em regular} in
the sense described above, i.e., if $B$ is any ball in ${\bf R}^s$ of
radius $t$, then $g^{-1}(B)$ can be covered by a bounded number of
balls of radius $t^\frac{s}{r}$.  This can be derived from the fact
that the $f_\epsilon$'s satisfy the same condition with uniform
bounds.

	The regularity condition is a kind of uniform and
scale-invariant way to say that $g$ has only bounded multiplicities.
It does not imply injectivity in general.  So what does the
injectivity of the mapping $f$ from Theorem \ref{stong's theorem}
really mean?  What does it imply at the continuous level?

	From regularity we also get that $g$ does not distort measures
very much.  That is, there is a constant $M$ so that
\begin{equation}
\label{distortion of measure}
	M^{-1} \, |E| \le |g(E)| \le M \, |E|
\end{equation}
for all, say, {\em compact} subsets of ${\bf R}^r$.  Here $|E|$ and
$|g(E)|$ denote the Lebesgue measures of the subsets $E$ and $g(E)$ of
${\bf R}^r$ and ${\bf R}^s$, respectively.  We do not really care
about compactness of $E$ here, it is just a convenient way to avoid
measurability problems.  One could as well work with outer measures
and have a similar estimate for {\em all} subsets $E$ of ${\bf R}^r$.
The upper bound in (\ref{distortion of measure}) follows from the fact
that $g$ is H\"older continous, while the lower bound can be derived
from the regularity property.  (For this it is convenient to use the
description of Lebesgue measure as Hausdorff measure of the correct
dimension.  Indeed, (\ref{distortion of measure}) has a more general
version for mappings between metric spaces and their effect on
Hausdorff measures.  This is spelled out in Lemma 12.3 of \cite{DS6}.)

	These bounds on the distortion of measure fit well with the
``measure-preserving'' property of bijections between discrete sets,
but they are less precise.  They permit bounded multiplicities, like
foldings, for instance.  Consider instead the condition that
\begin{equation}
\label{measure-preserving}
	|h(E)| = |E|
\end{equation}
for all compact subsets $E$ of ${\bf R}^r$, where $h : {\bf R}^r \to
{\bf R}^s$?  This is much stronger, closer to injectivity.  It implies
that disjoint compact sets in the domain are mapped to compact sets
which intersect in a set of measure zero.

	One can actually obtain mappings $h : {\bf R}^r \to {\bf R}^s$
which satisfy (\ref{measure-preserving}) from limits of rescalings of
$g$.  This follows from Proposition 12.42 of \cite{DS6}.  These new
mappings will also be H\"older continuous of the correct exponent and
satisfy the regularity property discussed above.  This
measure-preserving property (\ref{measure-preserving}) is much closer
to injectivity -- see
\cite{DS6} for more information about it -- but still it is not as
strong as the bijectivity in Theorem \ref{stong's theorem}.

	It seems that no matter what we do in the continuous case we
can never quite recapture fully the content of Theorem \ref{stong's
theorem}.  Infinite processes are involved in both cases, but the
condition of bijectivity in Theorem \ref{stong's theorem} does not
appear to have a reasonable formulation in the infinite case without
losing information, despite the fact that we can get a lot of structure
for the mappings in the continuous case.

	In neither case is the precise nature of the infinite processes
required so clearly understood.  One has examples and simple observations
about what kind of behavior is not possible, but deeper results are missing.

\section{A topological result about infinite complexity}

	By a {\em $k$-cell} we mean a topological space which is
homeomorphic to the Cartesian product of $k$ copies of the unit
interval $[0,1]$.

\begintheorem [\cite{Bry}]
\label{bryant's theorem}
Let $D$ be a compact subset of some ${\bf R}^n$ which is homeomorphic
to a cell of dimension $k$.  Let ${\bf R}^n / D$ denote the
topological space defined as a quotient of ${\bf R}^n$ by collapsing
$D$ to a point but leaving untouched the remaining elements of ${\bf
R}^n$.  Then $({\bf R}^n / D) \times {\bf R}$ is homeomorphic to
${\bf R}^{n+1}$.
\end{theorem}

	Thus one might say that shrinking a cell to a point does not
alter the ambient topology too painfully.  If the cell is a standard
one (a flat cube), then one can check directly that the quotient space
${\bf R}^n / D$ is homeomorphic to ${\bf R}^n$.  There are however
``wild'' cells for which this is not the case.  Theorem \ref{bryant's
theorem} is a bit remarkable for saying that after a single
stabilization the problem always goes away.

	This is a theorem about infinite complexity.  The assumptions
on $D$ allow infinite processes, and so does the conclusion.  It is a
positive statement about controlling infinite complexity.

	There are many other results of this general nature.  See
\cite{Dv2} for more information.

\section{Almost flat hypersurfaces}
\label{Almost flat hypersurfaces}

	Let $M$ be an hypersurface in ${\bf R}^{n+1}$.  For simplicity
we make the strong a priori assumptions that $M$ is a smooth embedded
submanifold which behaves well at infinity, e.g., is asymptotic to a
hyperplane.  We shall be interested in estimates which do not depend
on these assumptions in a quantitative way.

	Let $\epsilon > 0$ be given.  We shall say that $M$ is
$\epsilon$-flat \index{$\epsilon$-flat hypersurfaces} \index{almost
flat hypersurfaces} if it satisfies the following two properties.
The first asks that for each $x$ and $y$ on $M$ we can find a 
rectifiable curve $\gamma$ in $M$ which connects them and which
has length at most $(1+\epsilon) \, |x-y|$.  The second asks
that if $x \in M$ and $r > 0$ are given, then
\begin{equation}
	\biggl| \frac{|B(x,r) \cap M|}{\nu_n \, r^n} - 1 \biggr|
			\le \epsilon.
\end{equation}
Here $|B(x,r) \cap M|$ means the usual $n$-dimensional volume of
the intersection of the ball $B(x,r)$ with $M$, while $\nu_n$ denotes
the volume of the unit ball in ${\bf R}^n$.

	If $\epsilon = 0$ then $M$ must simply be an ordinary hyperplane
in ${\bf R}^{n+1}$.  If $\epsilon$ is small then these conditions imply
that $M$ is somehow pretty flat.

\beginquestion
Fix a positive integer $n$.  Is it true that for every $\eta > 0$
there is an $\epsilon > 0$ so that if $M$ is an $\epsilon$-flat
hypersurface in ${\bf R}^{n+1}$, then there is a bijection $h : {\bf R}^n
\to M$ which is $(1+\eta)$-bilipschitz, i.e.,
\begin{equation}
	(1+\eta)^{-1} \, |x-y| \le |h(x) - h(y)| \le (1+\eta) \, |x-y|
\end{equation}
for all $x, y \in {\bf R}^n$?
\end{question}

	It is easy to show that the converse to this is true.  More
precisely, if $M$ does admit such a $(1+\eta)$-bilipschitz
parameterization by ${\bf R}^n$, then $M$ is $\epsilon$-flat, where
$\epsilon \to 0$ (uniformly) as $\eta \to 0$.

	The a priori smoothness conditions here are convenient but not
really the point.  This question is really about the comparison
between two different measurements of infinite complexity.  Both
$\epsilon$-flatness and the $(1+\eta)$-bilipschitz allow infinite
processes, infinite amounts of spiralling, for instance.  They allow
one to have bumps of small but definite size at infinitely many scales
and locations at once.  But do they define practically the same
condition, as in the question above?

	This is not known.  There are some partial results, about
controlling the geometry of $\epsilon$-flat hypersurfaces, and finding
parameterizations which behave pretty well if not bilipschitzly, and
about analysis of functions and operators on such surfaces.  A general
theme is that practically anything that one might know for
hypersurfaces which admit $(1+\eta)$-bilipschitz parameterizations can
be proved directly for $\epsilon$-flat hypersurfaces, as long as one
allows some movement in the parameters, and one asks questions in a
reasonably simple language.  (One should not simply reformulate the
existence of a bilipschitz parameterization in other terms.)

	In affect one is able to obtain a lot of results about infinite
complexity in this case, even if the existence of good parameterizations
is not known.

	See \cite{casscI, casscII, pams, Sebil, Seapp, Sefred} for
more information about these and related questions.

\section{Big pieces of bilipschitz mappings}
\label{Big pieces of bilipschitz mappings}

\begintheorem
\label{guy's theorem}
Let $Q$ denote the unit cube in ${\bf R}^n$, and suppose that $f : Q \to
{\bf R}^n$ is Lipschitz with some constant $L$, i.e.,
\begin{equation}
	|f(x) - f(y)| \le L \, |x-y|
\end{equation}
for all $x, y \in Q$.  Assume also that the image of $f$ is of
substantial size, i.e., $|f(Q)| \ge \delta$, where $|f(Q)|$ denotes
the Lebesgue measure of $f(Q)$ and $\delta$ is some positive number.
Then there is an $\epsilon > 0$, depending only on $n, L$, and $\delta$,
and a closed subset $E$ of $Q$, such that $f$ satisfies is actually
bilipschitz on $E$, with
\begin{equation}
	|f(x) - f(y)| \ge \epsilon \, |x-y|
\end{equation}
for all $x, y \in E$.
\end{theorem}

	Of course the existence of a substantial bilipschitz piece
like this automatically ensures the existence of a lower bound on
$|f(Q)|$.  The fact that the converse holds, with uniform bounds,
is quite remarkable.

	This theorem was first proved in \cite{D4} (see also
\cite{D5}).  Another proof was given in \cite{J2}.  The latter proof was
simpler and gave a stronger result, to the effect that the domain $Q$
of $f$ could be decomposed into a bounded number of pieces, in such a
way that the restriction of $f$ to all but one of these pieces is
bilipschitz with a uniform bound, while the image of the remaining
piece under $f$ has small Lebesgue measure.  The more involved methods
of \cite{D4, D5} enjoy extra flexibility which enable them to be
applied in other situations.

	There are classical theorems in a somewhat similar spirit.
One can decompose $Q$ into a countable number of measurable subsets,
on each of which $f$ is either bilipschitz or has image of measure
$0$.  See \cite{Fe} for some results of this type, and for related
facts concerning the ``area'' of Lipschitz mappings.  These results
permit one to conclude that the jacobian of $f$ is bounded from below
on a set of substantial measure.  This is a weaker conclusion than the
existence of a substantial bilipschitz piece, as in the theorem above.
Indeed, one of the nice features of the theorem is that it lends
itself well to discretization, and indeed the discrete versions contain
essentially the same information as the continuous version.

	Thus Theorem \ref{guy's theorem} and the more classical
results which preceded it illustrate well the difference between
quantitative bounds and the more qualitative considerations of
countable unions and sets of measure zero.  The latter live within a
more symmetric language of infinite processes which is easier to use.
In Theorem \ref{guy's theorem} we must break some of the symmetry of
the language, but we then obtain a more concrete result which is
combinatorial in essence.  This fits well with the discussion of
Chapter \ref{About symmetry in language}.  For instance, a key point
in the classical theory is that an arbitrary countable union of sets
of measure zero also has measure zero.  It is not at all clear what
this means concretely, and one can certainly not say that arbitrary
finite unions of sets of small measure also have small measure.  To
get quantitative bounds one must go much further inside the structure
of the situation.

	Additional progress in finding large bilipschitz pieces of
Lipschitz mappings has been obtained in \cite{JKV}.

	The story of large bilipschitz pieces should be compared with
the one about $\epsilon$-flat hypersurfaces and the question about
bilipschitz parameterizations.  It should obviously be easier to find
bilipschitz pieces than global parameterizations, but that this is
true in a practical way is less obvious at first glance.  There is a
lot of technology available now for finding large bilipschitz pieces
of given spaces (additional references include \cite{DS2, DS4}), and
there is also now a reasonable collection of counterexamples
concerning the existence of global bilipschitz and other controlled
parameterizations \cite{Sebil, Seqs}.

	One can ask similar questions in a broader geometric context,
beyond spaces which have roughly Euclidean structure.  See \cite{DS6}
for more information.

\backmatter

\newpage
\addcontentsline{toc}{chapter}{Index}
\printindex

\end{document}